\documentclass[12pt,a4paper]{smfart}

\usepackage[all]{xy}
\usepackage[T1]{fontenc}
\usepackage[frenchb]{babel}
\usepackage{amsthm}
\usepackage{amsmath,amscd,amssymb,verbatim}

\usepackage{amsfonts}
\usepackage{amscd}

\swapnumbers

\def\limproj{\mathop{\oalign{$\mathrm{lim}$\cr
\hidewidth$\longleftarrow$\hidewidth\cr}}}

\newdimen\sumdim
\def\diagram#1\enddiagram{\vcenter{\offinterlineskip
\def\tvi{\vrule height 10pt depth 10pt width 10pt}\def\tvi{}
\halign{&\tvi\kern 5pt\hfil$\displaystyle ##$\hfil\kern 5pt
\crcr#1\crcr}}}
\def\hfl[#1][#2][#3]#4#5{\kern -#1
\sumdim=#2\advance\sumdim by #1\advance\sumdim by #3
\smash{\mathop{\hbox to\sumdim{\rightarrowfill}}
\limits^{\scriptstyle#4}_{\scriptstyle#5}}
\kern-#3}
\def\vfl[#1][#2][#3]#4#5%
{\sumdim=#1 \advance\sumdim by #2 \advance\sumdim by #3
\setbox1=\hbox{$\left\downarrow\vbox to .5\sumdim{}\right.$}
\setbox1=\hbox{\llap{$\scriptstyle #4$}\box1
\rlap{$\scriptstyle #5$}}
\vcenter{\kern -#1\box1\kern -#3}}
\def\Hfl #1#2{\hfl[0mm][12mm][0mm]{#1}{#2}}
\def\Vfl #1#2{\vfl[0mm][12mm][0mm]{#1}{#2}}

\def\Spm{\mathop{\mathrm{Spm}}}
\def\Spec{\mathop{\mathrm{Spec}}}
\def\Hom{{\mathrm{Hom}}}
\def\G{{\mathrm{Gr}}}
\def\Fitt{{\mathrm{Fitt}}}
\def\Spf{\mathop{\mathrm{Spf}}}
\def\dim{\mathop{\mathrm{dim}}}
\def\max{\mathop{\mathrm{max}}}

\def\mod{\mathop{\mathrm{modulo}}}

\def\ord{{\mathrm{ord}}}

\author{Julien Sebag}

\address{École Normale Supérieure, Département de mathématiques et applications, 
45 rue d'Ulm, 
75230 Paris cedex 05, France 
(UMR 8553 du CNRS)}
\email{Julien.Sebag@ens.fr}

\title{Intégration motivique sur les schémas formels}

\begin{document}

\maketitle

\let\!=\mathcal

\newtheorem{thm}[subsubsection]{Th\'eor\`eme}
\newtheorem{lem}[subsubsection]{Lemme}
\newtheorem{prop}[subsubsection]{Proposition}
\newtheorem{cor}[subsubsection]{Corollaire}
\newtheorem{defi}[subsubsection]{D\'efinition}
\newtheorem{exe}[subsubsection]{Exemple}
\newtheorem{rap}[subsubsection]{Rappel}
\newtheorem{rem}[subsubsection]{Remarque}
\newtheorem{lemcle}[subsubsection]{Lemme clé}

\newtheorem{slem}[subsubsection]{Sous-lemme}

\newtheorem{nota}[subsubsection]{Notations}
\newtheorem{exes}[subsubsection]{Exemples}

\newdimen\sumdim
\def\diagram#1\enddiagram{\vcenter{\offinterlineskip
\def\tvi{\vrule height 10pt depth 10pt width 10pt}\def\tvi{}
\halign{&\tvi\kern 5pt\hfil$\displaystyle ##$\hfil\kern 5pt
\crcr#1\crcr}}}
\def\hfl[#1][#2][#3]#4#5{\kern -#1
\sumdim=#2\advance\sumdim by #1\advance\sumdim by #3
\smash{\mathop{\hbox to\sumdim{\rightarrowfill}}
\limits^{\scriptstyle#4}_{\scriptstyle#5}}
\kern-#3}
\def\vfl[#1][#2][#3]#4#5%
{\sumdim=#1 \advance\sumdim by #2 \advance\sumdim by #3
\setbox1=\hbox{$\left\downarrow\vbox to .5\sumdim{}\right.$}
\setbox1=\hbox{\llap{$\scriptstyle #4$}\box1
\rlap{$\scriptstyle #5$}}
\vcenter{\kern -#1\box1\kern -#3}}
\def\Hfl #1#2{\hfl[0mm][12mm][0mm]{#1}{#2}}
\def\Vfl #1#2{\vfl[0mm][12mm][0mm]{#1}{#2}}


\section{Introduction}

Dans leurs articles \cite{dl1} et \cite{dl2}, Denef et Loeser d\'eveloppent et étudient la th\'eorie de l'\textit{int\'egration motivique}, introduite par Kontsevich lors d'un s\'eminaire \`a Orsay (\textit{cf} \cite{K}). Cette nouvelle th\'eorie de l'int\'egration se r\'evèle \^etre un outil puissant dans l'étude de la g\'eom\'etrie birationnelle des variétés algébriques sur un corps $k$ de caractéristique 0. Rappelons rapidement les grandes étapes de la construction de ces intégrales et les idées sous-jacentes à quelques-uns des résultats obtenus par cette théorie. Si $X$ est une variété algébrique sur un corps $k$ de caractéristique 0, on lui associe, de manière fonctorielle, un pro-$k$-schéma, qui est encore un schéma (non de type fini en général), le \textit{schéma des arcs sur $X$}, noté $\!L(X)$. Sur ce pro-$k$-schéma, on définit une mesure $\mu_X$ à valeurs dans un anneau de motifs virtuels, que l'on note $\!M$. Cette mesure motivique est un analogue géométrique de la mesure $p$-adique sur les variétés différentielles $p$-adiques. Les intégrales définies à partir de cette mesure vérifient alors un théorème de changement de variables pour les $k$-morphismes de schémas $h:Y\rightarrow X$ propres et birationnels, qui permet de calculer les intégrales sur $\!L(X)$ en fonction d'intégrales sur $\!L(Y)$. C'est essentiellement par ce principe que l'on peut construire de nouveaux invariants algébriques (à valeurs dans cet anneau de motifs virtuels, ou plus exactement dans le séparé complété de cet anneau pour une filtration) et que l'on peut, par exemple, (re)-démontrer le théorème de Batyrev qui affirme que deux variétés de Calabi-Yau ont même nombre de Hodge et même structure de Hodge.

Dans cet article, nous généralisons cette théorie de l'intégration motivique aux sch\'emas formels \textit{sttf} sur le spectre formel d'un anneau de valuation discr\`ete complet $R$, de corps résiduel $k$ parfait. En particulier, le corps $k$ peut être de caractéristique positive. L'absence du morphisme d'inclusion $k\hookrightarrow R$, dans le cas où $R$ est un anneau d'inégale caractéristique, nous conduit à définir un analogue du schéma des arcs d'une variété $X$. Si $\!X$ est un $R$-schéma formel \textit{sttf}, le \textit{schéma de Greenberg} $\G(\!X)$ joue ce rôle. Il est important de noter que, dans le cas où $R=k[[t]]$, avec $k$ de caractéristique 0, et $\!X$ est le complété $\pi$-adique d'une variété algébrique sur $k$, les deux $k$-schémas $\G(\!X)$ et $\!{L}(X)$ sont canoniquement isomorphes. En outre, pour tout $R$-schéma formel \textit{sttf}, si $F$ est une extension parfaite de $k$, les $F$-points de $\G(\!X)$ s'interprètent naturellement comme des points de la fibre générique $X_K$ de $\!X$. Comme dans le cas du schéma des arcs, nous définissons une mesure $\mu_\!X$ à valeurs dans le même complété de cet anneau de motifs virtuels et l'anneau booléen des \textit{ensembles mesurables}, qui sont des ``approximations'' par certaines parties constructibles de $\G(\!X)$ élémentaires, que l'on appelle \textit{cylindres}. De manière analogue aux théories classiques d'intégration, les ensembles mesurables possèdent des propriétés de stabilité par image directe et inverse sous certains $R$-morphismes $h:\!Y\rightarrow \!X$, que l'on nomme \textit{tempérés}, et qui correspondent aux morphismes propres et birationnels du cas algébrique. Nous construisons l'intégrale motivique des fonctions intégrables. Parmi celles-ci, les fonctions \textit{exponentiellement intégrables} permettent d'exprimer, de manière naturelle, des phénomènes géométriques comme l'appartenance à un sous-$R$-schéma formel fermé, la lissité d'un $R$-morphisme de schémas formels $\ldots$ et jouent donc un rôle central dans cette théorie. Enfin, nous démontrons deux théorèmes de changement de variables du même type que celui énoncé dans \cite{dl1} ou \cite{dl2}.

Les principaux résultats de ce travail sont utilisés de manière fondamentale dans \cite{ls}, qui, en se pla\c{c}ant du point de vue rigide, déduit, des constructions et théorèmes de cet article, certaines applications en géométrie birationnelle des dégénérescences des variétés algébriques et rigides. En particulier, la théorie que l'on développe ici apparaît, au regard des résultats de \cite{ls}, comme une généralisation de l'intégration motivique usuelle et de l'intégration $p$-adique. Deux exemples le confirment. Soit $X_K$ la fibre générique d'un $R$-schéma formel $\!X$ \textit{sttf}.

\begin{enumerate}

\item Si $\omega$ est une forme jauge sur $X_K$, son intégrale converge déjà dans $\!M$. En outre, nous démontrons que sa classe dans $\!M/(\mathbf{L}-1)\!M$ ne dépend que de $X_K$ et non de la forme jauge choisie pour calculer cette intégrale. Cet élément de $\!M$, noté $\lambda(X_K)$, est égal à la classe de la fibre spéciale de tout  modèle de Néron faible de $X_K$. Il faut également remarquer que, quand $K$ est une extension finie de $\mathbf{Q}_p$, $\lambda(X_K)$ se spécialise en l'invariant de Serre \cite{s} évalué sur la variété analytique $p$-adique sous-jacente à $X_K$ (\textit{cf} corollaire 4.6.3 de \cite{ls}).

\item Si $X_K$ est l'espace analytique associé une variété de Calabi-Yau sur $K$, et si $X_K$ admet un $R$-modèle propre et lisse, l'intégrale calculée pour un générateur de $\Omega_{X_K}^d$ est égale à la classe de la fibre spéciale de ce modèle dans $\!M$. En particulier, les classes des fibres spéciales de deux tels modèles co\"\i ncident dans $\!M$, ce qui peut être interprété comme l'analogue du résultat de Batyrev pour les variétés de Calabi-Yau (\textit{cf} \cite{ba0}).
\end{enumerate}

Cet article est organisé de la manière suivante: le chapitre 2 rappelle les bases de la géométrie formelle. Le foncteur de Greenberg est construit et étudié au chapitre 3. Les chapitres 4, 5, et 6 développent, de manière assez complète, les notions et les propriétés de la mesure motivique, des cylindres et des ensembles mesurables. Le chapitre 7 est exclusivement consacré à la construction de l'intégrale motivique, à l'énoncé et à la démonstration des deux théorèmes de changement de variables.

Nous remercions Fran\c{c}ois Loeser de nous avoir proposé ce sujet et de l'aide qu'il nous a apport\'ee lors des discussions que nous avons eues ensemble. Nous souhaitons également remercier Antoine Chambert-Loir de son aide lors de la rédaction de cet article et de ses commentaires sur une première version de ce texte. Plus largement, nous les remercions pour l'attention qu'ils nous ont manifesté depuis quelques années déjà.


\section{Préliminaires}



\subsection{Les notations.}


Dans cet article, l'anneau $R$ est un anneau de valuation discrète complet de corps résiduel $k$ et de corps des fractions $K$. Le corps $k$ est supposé parfait. On fixe $\pi$ une uniformisante de cet anneau. On note $R_n:=R/(\pi)^{n+1}$ pour tout $n\in\mathbf{N}$. On note $\mathbb{D}:=\Spf R$ le spectre formel de l'anneau $R$. L'espace topologique sous-jacent à cet espace localement annelé est constitué d'un unique point, qui s'identifie au point fermé du schéma $\Spec R$ (\textit{cf} \cite{ega} \S 10). On note $\mathbb{B}^{N}_R:=\Spf R\{x_1,\ldots,x_N\}$ l'espace affine formel de dimension $N$ sur $R$.


\subsection{La description des anneaux de valuation discrète complets.}
\label{description de R}


Tout élément $x\in R$ s'écrit de manière unique sous la forme:

\begin{enumerate}

\item $x=\sum_{i=0}^\infty a_i t^i$, avec $a_i\in k$ et $t:=\pi$ dans le cas d'égale caractéristique.

\item $x=\sum_{i=0}^\infty\tau(a_i^{p^{-i}})p^i$, avec $a_i\in k$, $\tau$ le morphisme de Teichm\"{u}ller et $p=\pi$ dans le cas d'inégale caractéristique, non-ramifié.

\item $x=\sum_{i=0}^\infty\tau(a_i^{p^{-q(i)}})p^{q(i)}\pi^{r(i)}$, avec $a_i\in k$, $\tau$ le morphisme de Teichm\"{u}ller et les entiers $q(i)$ et $r(i)$ donnés par la division euclidienne de $i$ par $e$, $i=q(i)e+r(i)$, avec $e$ l'indice de ramification, dans le cas d'inégale caractéristique ramifié.

\end{enumerate}


\subsection{Les premières notions de Géométrie formelle.}


\subsubsection{La définition des schémas formels}

Un $R$-schéma formel désignera toujours un $R$-schéma formel topologiquement de type fini, ce que l'on notera \textit{ttf}. Parfois l'hypothèse de séparation sera nécessaire et nos $R$-schémas formels \textit{ttf} seront séparés, ce que l'on notera \textit{sttf} (\textit{cf} \cite{ega} \S 10). On notera $\underline{Form}^{\mathrm{ttf}}_{/R}$ et $\underline{Form}^{\mathrm{sttf}}_{/R}$ les catégories correspondantes.

Un objet de $\underline{Form}^{\mathrm{ttf}}_{/R}$ est un espace localement annelé $(\!X,\!O_{\!X})$ en $R$-algèbres topologiques, qui induit la donnée, pour tout $n\geq 0$, d'un $R_n$-sch\'ema $X_n=(\!X,\!O_\!X\otimes_R R_n)$. Le $k$-schéma $X_0$ est appelé fibre spéciale du $R$-schéma formel $\!X$. En tant qu'espaces topologiques, $\!X$ et $X_0$ sont isomorphes et $\!O_\!X:=\limproj\!O_{X_n}$. On a $X_n=X_{n+1}\otimes_{R_{n+1}} R_n$ et $\!X$ est canoniquement isomorphe à la limite inductive des schémas $X_n$ dans la catégorie des espaces localement annelés. En outre, un tel objet est localement isomorphe à un $R$-schéma formel affine $\Spf A$, où $A$ est une $R$-algèbre $\pi$-adique, topologiquement isomorphe à un quotient de l'anneau des séries formelles restreintes $R\{T_1,\ldots,T_N\}$.

Si $\!X$ et $\!Y$ sont deux $R$-schémas formels \textit{ttf}, on note $\Hom_R(\!Y, \!X)$ l'ensemble des $R$-morphismes de schémas formels:
$$
\xymatrix{\!Y\ar[rr]\ar[dr]&&\!X\ar[dl]\\
&\mathbb{D}&}
$$
\textit{i.e.} l'ensemble des morphismes entre les $R$-espaces localement annelés sous-jacents. Autrement dit, $\underline{Form}^{\mathrm{ttf}}_{/R}$ est une sous-catégorie pleine de la catégorie des espaces localement annelés sur $R$. Localement de tels morphismes sont simplement des $R$-morphismes continus d'algèbres topologiques entre les anneaux des sections globales. On peut montrer (\textit{cf} proposition 10.6.9 de \cite{ega}) que l'application canonique $\Hom_R(\!Y,\!X)\rightarrow\limproj\Hom_{R_n}(Y_n,X_n)$ est une bijection.

\begin{rem}
Soit 
$$
h:\!Y=\Spf R\{x_1,\ldots,x_m\}/\!I\rightarrow\!X=\Spf R\{x_1,\ldots,x_n\}/\!J
$$ 
un $R$-morphisme de schémas formels. Dans ce cas, la donnée de $h$ est équivalente à celle d'un $R$-morphisme d'algèbres. En effet, la $R$-linéarité impose que le morphisme entre les anneaux des sections globales soit continu pour les topologies $\pi$-adiques considérées.
\end{rem}

La catégorie $\underline{Form}^{\mathrm{sttf}}_{/R}$ est simplement la sous-catégorie pleine de $\underline{Form}^{\mathrm{ttf}}_{/R}$, dont les objets sont séparés (\textit{cf} \cite{ega} \S 10). Un $R$-schéma formel \textit{sttf} $\!X$ est dit \textit{admissible} s'il est plat sur $R$ (ce qui est équivalent au fait que le faisceau $\!O_\!X$ soit sans $\pi$-torsion). On notera $\underline{Form}_{/R}^{\textrm{Adm}}$ cette catégorie, qui est une sous-catégorie pleine de $\underline{Form}^{\mathrm{sttf}}_{/R}$.

\begin{lem}
\label{neron=sttf}
Soit $\!X$ un $R$-schéma formel ttf et $j:\!U\hookrightarrow\!X$ une immersion ouverte. Alors, si $\!X$ est séparé et $\!U$ quasi-compact, le $R$-morphisme $j$ est rétro-compact, \textit{i.e.} pour tout ouvert affine $\!V$ de $\!X$, $\!V\cap \!U$ est quasi-compact. En particulier, le $R$-schéma formel $\!U$ est sttf.
\end{lem}

\begin{proof}
Soit $(\!U_i)$ un recouvrement ouvert fini de $\!U$ par des sous-$R$-schémas formels affines \textit{ttf}. Comme $\!U$ est ouvert et $\!X$ séparé, $\!U$ est également séparé. Le schéma formel $\!V\cap \!U$ est quasi-compact, comme réunion finie d'ouverts quasi-compacts.
\end{proof}

On dit qu'un $R$-schéma formel est \textit{réduit} si sa fibre spéciale est un schéma réduit. Si $\!X$ est un $R$-schéma formel, on notera $\!X_{\mathrm{red}}$ l'unique sous-schéma formel ferm\'e et réduit de $\!X$.

\subsubsection{La notion de fibre générique d'un schéma formel}

Soit $\!X$ un $R$-schéma formel \textit{ttf} et $\!Z\hookrightarrow\!X$ un sous-schéma fermé de $\!X$, défini par le sous-faisceau $\!A\subset\!O_\!X$. On appelle \textit{éclatement admissible sur $\!X$ de centre $\!Z$} la donnée d'un $R$-schéma formel $\!X'$ et d'un $R$-morphisme de schémas formels $\sigma:\!X'\rightarrow\!X$ vérifiant la propriété universelle suivante: si $\psi:\!Y\rightarrow\!X$ est un $R$-morphisme de schémas formels tel que $\!A.\!O_\!Y$ soit inversible sur $\!Y$, alors il existe un $R$-morphisme $\psi':\!Y\rightarrow\!X'$ tel que $\psi=\sigma\circ\psi'$. Si $\!X$ est \textit{sttf} (\textit{resp. admissible}),  le $R$-schéma formel $\!X'$ l'est aussi (\textit{cf} \cite{bl1}).

La localisation de la catégorie $\underline{Form}^{\mathrm{ttf}}_{/R}$ (\textit{resp.} $\underline{Form}^{\mathrm{sttf}}_{/R}$) par rapport aux éclatements admissibles est équivalente à la catégorie des $K$-espaces rigides de type fini et quasi-séparés (\textit{resp.} séparés) au sens de Kiehl, que l'on notera $\underline{Rig}^{\mathrm{qsqc}}_{/K}$ (resp. $\underline{Rig}^{\mathrm{sqc}}_{/K}$) (\textit{cf} \cite{r} et \cite{bl1}). On peut remarquer que la localisée de $\underline{Form}^{\mathrm{Adm}}_{/R}$ par rapport aux éclatements admissibles est équivalente à $\underline{Rig}^{\mathrm{sqc}}_{/K}$.

Par analogie avec le cas de schémas usuels, on appellera fibre générique l'image d'un $R$-schéma formel \textit{ttf} par le foncteur de localisation. Ce foncteur sera noté $\mathrm{rig}$ et si $\!X\in\underline{Form}^{\mathrm{ttf}}_{/R}$ on notera $\!X_{\mathrm{rig}}$ son image par le foncteur $\mathrm{rig}$ (ou parfois $X_K$). De même, si $f:\!Y\rightarrow\!X$ est un $R$-morphisme de schémas formels, on notera $f_{\mathrm{rig}}$ son image par le foncteur $\mathrm{rig}$.

\subsubsection{La notion de dimension d'un schéma formel}
\label{notion de lissite}

Soit $\!X$ un $R$-schéma formel. On appelle \textit{dimension} de $\!X$ l'entier $\dim \!X$ défini comme la dimension de la fibre spéciale $X_0$ de $\!X$. En particulier, si $\!X:=\Spf A$ est affine, on a la relation
$$
\dim \!X=\textrm{dim}_{\textrm{K}} A -1
$$
où $\dim_{\textrm{K}} A$ est la dimension de Krull de l'anneau $A$.

Soit $X_K$ un $K$-espace rigide quasi-séparé et quasi-compact. On appelle \textit{dimension} de $X_K$ et l'on note $\dim X_K$ l'entier naturel défini comme la borne supérieure des entiers $\dim_{\mathrm{K}} \!O_{X_K,x}$ pour $x\in X_K$.

Le lemme suivant, démontré par Oesterlé dans \cite{o1} (\textit{cf} lemme 1 \S 3), permet de relier ces deux notions.
\begin{lem}
\label{theoreme de Osterle}
Soit $R$ un anneau de valuation discrète complet, $k$ son corps résiduel et $K$ son corps des fractions. Notons $A$ l'anneau $R\{T_1,\ldots,T_N\}$ des séries formelles restreintes à coefficients dans $R$, $A_K=A\otimes_R K=K\{T_1,\ldots,T_N\}$ et $A_0=A\otimes_R k=k[T_1,\ldots,T_N]$. Soit $I_K$ un idéal de $A_K$. Notons $I:=I_K\cap A$ et $I_0:=I\otimes_R k$. Supposons que l'anneau $A_K/I_K$ soit équidimensionnel, de dimension de Krull $d$. Alors l'anneau $A_0/{I_0}$ est équidimensionnel, de dimension de Krull $d$. On a $I_K=A_K$ si et seulement si ${I_0}=A_0$.
\end{lem}

\begin{lem}
Si $X_K=\Spm A$ est un $K$-espace rigide affino\"{i}de, la dimension de $X_K$ est égale à la dimension de Krull de $A$.
\end{lem}

\begin{proof}
Ce lemme résulte de la proposition 7.3/8 de \cite{bgr} qui assure que l'anneau $\!O_{X_K,x}$ et le localisé de $A$ en $x$, $A_x$, ont même dimension de Krull.
\end{proof}

\begin{cor}
Soit $\!X$ un $R$-schéma formel sttf de fibre générique $X_K$. Si $X_K$ est équidimensionnelle de dimension $d$, alors $X_0$ est équidimensionnel de dimension $d$.
\end{cor}

\begin{proof}
On peut supposer que $\!X$ est affine. Dans ce cas le corollaire découle du lemme \ref{theoreme de Osterle}.
\end{proof}

Un $R$-schéma formel \textit{sttf} est dit \textit{équidimensionnel} de dimension $d$ (ou \textit{de pure dimension $d$}) si sa fibre générique $X_K$ est équidimensionnelle de dimension $d$.

\subsubsection{La notion de lissité pour les schémas formels}
\label{lissite}

Un $R$-morphisme de schémas formels $f:\!Y\rightarrow\!X$ est \textit{lisse au point $y\in Y_0$} de dimension relative $d$ si:

\begin{enumerate}

\item $f$ est plat en $x$.

\item Le $k$-morphisme induit, $f_0:Y_0\rightarrow Y_0$, est lisse en $y$ de dimension relative $d$, au sens usuel.

\end{enumerate}

Il revient au même (\textit{cf} lemme 1.2 de \cite{bl1}) de demander que pour tout $n\in\mathbf{N}$, le morphisme induit $f_n:Y_n\rightarrow X_n$ soit lisse de dimension relative $d$. On dit qu'un morphisme est lisse s'il est lisse en tout point $x\in X_0$. 

Un $R$-schéma formel $\!X$ est dit lisse (en $x\in X_0$) si le morphisme structural $\!X\rightarrow\Spf R$ est lisse (en $x\in X_0$) au sens précédent.

Soit $\!X$ un $R$-schéma formel \textit{ttf} de dimension $d$. On note $\!X_{\mathrm{sing}}$ l'unique sous-$R$-schéma formel fermé réduit défini par le radical du $d$-ème idéal de Fitting de $\Omega_{\!X/R}^1$ (\textit{cf} \cite{bl2} \S 3). En particulier, un $R$-schéma formel \textit{ttf} plat est lisse en $x\in X_0$ (\textit{resp.} lisse) si et seulement si $x$ n'appartient pas à la fibre spéciale de $\!X_{\mathrm{sing}}$ (\textit{resp.} $\!X_{\mathrm{sing}}=\emptyset$).

\begin{lem}
\label{autre definition du lieu singulier}
Soit $\!X:=\Spf A\hookrightarrow\mathbb{B}^N_R$ un $R$-schéma formel affine admissible, avec $A:=R\{x_1,\ldots,x_N\}/\!I$. Si $f_1,\ldots,f_m\in \!I$ et $m\leq N$, soit $\Delta(f_1,\ldots,f_m)$ l'idéal de $R\{x_1,\ldots,x_N\}/\!I$ engendré par les $m\times m$-mineurs de la matrice jacobienne $(\partial f_j/\partial x_i)$. Soit $H_{A/R}$ le radical dans $A$ de l'idéal $\sum\Delta(f_1,\ldots,f_m)((f_1,\ldots,f_m):\!I)$, où la somme est prise sur l'ensemble des ensembles finis $(f_1,\ldots,f_m)\in\!I$. Alors $\!X$ est lisse en $x\in X_0$ si et seulement si $x$ n'appartient pas au sous-$R$-schéma formel fermé de $\!X$ défini par l'idéal $H_{A/R}$.
\end{lem}

\begin{proof}
Comme $\!X$ est plat, la lissité se lit sur la fibre spéciale $X_0$, ainsi que l'appartenance au sous-$R$-schéma formel fermé. Le lemme est donc une conséquence du théorème 4.1 de \cite{pop}.  
\end{proof}

On dira qu'un $R$-schéma formel $\!X$ \textit{ttf} est \textit{génériquement lisse} si sa fibre générique $X_K$ est lisse comme $K$-espace rigide, \textit{i.e.} vérifie un critère jacobien identique à celui vérifié par les schémas lisses (\textit{cf} \cite{bl3}).


\section{Foncteur de Greenberg et (pro)-schéma associé à un schéma formel}



\subsection{La construction}


Dans ce paragraphe, nous allons construire un foncteur
$$
\G:\underline{Form}^{\mathrm{sttf}}_{/R}\rightarrow\underline{ProSch}_{/k}
$$
de la catégorie des $R$-schémas formels \textit{sttf} dans la catégorie des pro-$k$-schémas.

Rappelons quelques résultats de \cite{g1} et \cite{blr} \S 9.6. Pour tout $n\geq 0$, considérons le foncteur $h^\ast_n$ qui à un $k$-schéma associe

\begin{enumerate}

\item le $k$-schéma $h_n^\ast(T):=T\times_k R_n$, en égale caractéristique.

\item l'espace localement annelé $h_n^\ast(T)$, dont l'espace topologique sous-jacent est $T$ et le faisceau structural est $\!Hom(T,\!R_n)$, où $\!R_n$ est le $k$-schéma en anneaux dont l'anneau des points rationnels est $R_n$ (\textit{cf} \cite{blr} et \cite{ls}), dans les cas d'inégale caractéristique. Dans le cas où $R=W(k)$, ce $k$-schéma $\!R_n$ est simplement le schéma en anneaux des vecteurs de Witt de longueur $n$.

\end{enumerate}

Soit $A$ une $k$-algèbre. On pose $L(A)=A$ si $R$ est un anneau d'égale caractéristique et $L(A)=W(A)$ si $R$ est un anneau d'inégale caractéristique.

En particulier, pour toute $k$-algèbre $A$,
$$
h_n^\ast(\Spec A)=:h_n^\ast(A)=\Spec(R_n\otimes_{L(k)} L(A)).
$$

Dans l'article \cite{g1}, Greenberg a montré que, pour tout $R_n$-schéma $X_n$ localement de type fini, le foncteur
$$
T\mapsto\Hom_R(h_n^\ast(T), X_n)
$$
de la catégorie des $k$-schémas dans celle des ensembles est représentable par un $k$-schéma $\G_n(X_n)$ qui est localement de type fini (\textit{i.e.} le foncteur $h^\ast_n$ admet un adjoint à droite que l'on note $\G_n$) et vérifie, pour toute $k$-algèbre $A$,
$$
\G_n(X_n)(A)\simeq X_n(R_n\otimes_{L(k)} L(A)).
$$
En particulier, si $A=k$, on a 
$$
\G_n(X_n)(k)\simeq X_n(R_n).
$$

Soit $\!X$ un $R$-schéma formel \textit{sttf}. Les morphismes canoniques $R_{n+1}\rightarrow R_n$, induisent alors, notamment par adjonction, pour tout $n\in\mathbf{N}$, des $k$-morphismes
$$
\theta_{n}^{n+1}:\G_{n+1}(X_{n+1})\rightarrow\G_n(X_n)
$$
qui font de la suite $(\G_n(X_n))_{n\in\mathbf{N}}$ un système projectif dans la catégorie des $k$-schémas de type fini et séparés (car, pour tout $n$, le $k$-schéma $X_n$ est séparé). Enfin, comme un recouvrement (fini) de $\!X$ par des ouverts formels affines induit un recouvrement des $\G_n(X_n)$ (\textit{cf} proposition \ref{recollement}/2), ces morphismes de transition sont affines et leur limite $\limproj\G_n(X_n)$ existe dans la catégorie des $k$-schémas.

Soit $T$ un $k$-schéma. On définit $h^\ast(T)$ l'espace localement annelé, dont l'espace topologique sous-jacent est $T$ et le faisceau structural $\limproj\!Hom_k(T,\!R_n)$. Cet espace localement annelé appartient à la catégorie des espaces localement annelés sur $\mathbb{D}$.

\begin{lem}
\label{groth-Gr}
Pour tout $R$-schéma formel sttf $\!X$, l'application canonique
$$
\Hom_\mathbb{D}(h^\ast(T),\!X)\rightarrow\limproj\Hom_{R_n}(h^\ast_n(T),X_n)
$$
est une bijection.
\end{lem}

\begin{proof}
Découle des définitions.
\end{proof}

\begin{prop}
Soit $\!X$ un $R$-schéma formel sttf. Le foncteur
$$
T\mapsto\Hom_R(h^\ast(T),\!X)
$$
de la catégorie des $k$-schémas dans celle des ensembles est représentable par un $k$-schéma canoniquement isomorphe à $\limproj\G_n(X_n)$.
\end{prop}

\begin{proof}
Ceci découle de la propriété universelle des $\G_n(X_n)$ et du lemme \ref{groth-Gr}.
\end{proof}

On posera alors 
$$
\G_n(\!X):=(\G_n(X_n))_{\mathrm{red}}
$$ 
et 
$$
\G(\!X):=\limproj(\G_n(X_n)_{\mathrm{red}})\simeq(\limproj\G_n(X_n))_{\mathrm{red}}.
$$
On a donc construit un foncteur $\G$ de la catégorie des $R$-schémas formels \textit{sttf} dans la catégorie des pro-$k$-schémas, comme composé des foncteurs définis ci-dessus.

En particulier, pour toute extension $F$ de $k$, les applications canoniques
$$
\G_n(\!X)(F)\rightarrow\G_n(X_n)(F)
$$
et
$$
\G(\!X)(F)\rightarrow(\limproj\G_n(X_n))(F)
$$
sont des bijections.

Dans le cas d'égale caractéristique, $R\simeq k[[\pi]]$. Soit $X$ une variété algébrique sur $k$. On peut alors considérer pour tout $n\in\mathbf{N}$, le $R_n$-schéma $X\times_k R_n$. Il découle des définitions que $\G_n(X\times_k R_n)$ est isomorphe au $k$-schéma $\!L_n(X)$ considéré dans \cite{dl1}. On en déduit en particulier que $\!L(X)\simeq\G(X\hat\otimes_k R)$, où $X\hat\otimes_k R$ est simplement le complété $\pi$-adique de $X$.

\begin{exe}
Il découle des définitions que le $k$-schéma $\G_0(\!X)$ est canoniquement isomorphe à $(X_0)_{\mathrm{red}}$.
\end{exe}


\subsection{Les notations}


\begin{enumerate}

\item

Soit $A$ une $k$-algèbre. On pose $L(A)=A$ si $R$ est un anneau d'égale caractéristique et $L(A)=W(A)$ si $R$ est un anneau d'inégale caractéristique. On notera par $R_A$ l'anneau $R_A:=R\otimes_{L(k)}L(A)$. Si $F$ est un corps parfait contenant $k$, on notera $K_F$ le corps des fractions de $R_F$. On peut remarquer alors, que, comme $k$ est supposé parfait, l'extension $R\rightarrow R_F$ est d'indice de ramification 1 au sens de \cite{blr} \S 3.6. En particulier, on a une bijection canonique
$$
\G(\!X)(F)\simeq\!X(R_F).
$$

\item

Pour tout $n\in\mathbf{N}$, les morphismes canoniques seront toujours notés de la manière suivante:

$$
\xymatrix{\G(\!X)\ar[r]^{\pi_{n,\!X}}\ar[dr]_{\pi_{n-1,\!X}}& \G_n(\!X)\ar[d]^{\theta^n_{n-1}} \\
& \G_{n-1}(\!X).}
$$
Les morphismes $\pi_{n,\!X}$ (ou $\pi_n$) sont les morphismes de troncation. Les morphismes $\theta_n^{n-1}$ sont les morphismes de transition. 

\item

Soit $\!Y$ un $R$-schéma formel \textit{sttf}. Soit $h:\!Y\rightarrow\!X$ un $R$-morphisme de schémas formels. On notera encore $h$ le $k$-morphisme de schémas $\G(h)$, et on notera $h_n$ le $k$-morphisme $\G_n(h)$, pour tout $n\in\mathbf{N}$. Ces morphismes $h$ et $h_n$ rendent commutatif le diagramme suivant:
$$
\xymatrix{\G(\!Y)\ar[r]^{h}\ar[d]_{\pi_{n,\!Y}}&\G(\!X)\ar[d]^{\pi_{n,\!X}}\\
\G_n(\!Y)\ar[r]_{h_n}&\G_n(\!X)}
$$

\end{enumerate}


\subsection{Les propriétés et la définition par recollement}


\begin{prop}

\begin{enumerate}
\label{recollement}

\item Le foncteur $\G$ préserve les immersions ouvertes et fermées, les produits fibrés, et transforme $R$-schémas formels affines en $k$-schémas affines.

\item Soit $\!X$ un $R$-schéma formel sttf et soit $(\!O_i)_{i\in J}$ un recouvrement fini de $\!X$ par des sous-schémas formels ouverts affines. Il existe des isomorphismes canoniques $\G(\!O_i\cap\!O_j)\simeq\G(\!O_i)\cap\G(\!O_j)$ et le $k$-schéma $\G(\!X)$ est obtenu en recollant les $\G(\!O_i)$.
\end{enumerate}
\end{prop}

\begin{proof}

Les assertions \ref{recollement}/1 et \ref{recollement}/2 pour $\G$ se déduisent directement de celles pour $\G_n$. L'assertion \ref{recollement}/1 est démontrée dans \cite{g1}. Il nous suffit donc de prouver l'assertion \ref{recollement}/2 pour $\G_n$ et $n\in\mathbf{N}$. L'existence d'isomorphismes canoniques découle de la propriété universelle. L'assertion \ref{recollement}/2 s'ensuit.
\end{proof}

\begin{lem}
Soit $\!X$ un $R$-schéma formel sttf. Alors le $k$-morphisme canonique de schémas 
$$
\G(\!X_{\mathrm{red}})\hookrightarrow\G(\!X)
$$
est un isomorphisme.
\end{lem}

\begin{proof}
Soit $x\in\G(\!X)$. Par adjonction, il correspond, à un point géométrique au-dessus de $x$, un $\mathbb{D}$-morphisme de schémas formels $\varphi:\Spf R_{k'}\rightarrow\!X$, où $k'$ est une clôture algébrique de $\kappa(x)$, le corps résiduel de $x$ dans $\G(\!X)$. Comme $\Spf R_{k'}$ est réduit, le morphisme $\varphi$ se factorise par $\!X_{\mathrm{red}}$. En particulier, $x\in \G(\!X_{\mathrm{red}})$. Les deux $k$-schémas réduits de $\G(\!X)$ et $\G(\!X_{\mathrm{red}})$ ont donc le même espace topologique sous-jacent. L'assertion en découle.
\end{proof}

Rappelons la définition d'une algèbre de Tate analytiquement séparable (\textit{cf} \cite{sch}).

\begin{defi}
Soit $A$ une algèbre de Tate définie par un idéal 
$$
\mathfrak{a}\subset K\{X_1,\ldots,X_N\}.
$$
On dit que $A$ est analytiquement séparable si l'anneau
$$
\overline{K}\{X_1,\ldots,X_N\}/\mathfrak{a}\overline{K}\{X_1,\ldots,X_N\}
$$
est réduit, où $\overline{K}$ est le complété d'une clôture algébrique de $K$. Dire que $A$ n'est pas analytiquement séparable équivaut à dire que $K$ est de caractéristique $p>0$ et que $K^{1/p}\{X_1,\ldots,X_N\}/\mathfrak{a}K^{1/p}\{X_1,\ldots,X_N\}$ n'est pas réduit.
\end{defi}

\begin{prop}
\label{schappa}
Soit $\!X$ un $R$-schéma formel admissible et réduit de pure dimension $d$. Soit $\!X_{\mathrm{sing}}$ le sous-$R$-schéma formel fermé de $\!X$ défini en \S \ref{notion de lissite}. Soient $X_K$ la fibre générique de $\!X$ et $(\!X_{\mathrm{sing}})_K$ celle de $\!X_{\mathrm{sing}}$. Alors soit $(\!X_{\mathrm{sing}})_K$ est de codimension au moins 1 dans $X_K$; soit, pour tout ouvert affine $\!U\rightarrow\!X$ sttf, il existe un sous-$R$-schéma formel $\!Z_\!U\hookrightarrow\!U$ tel que:

\begin{enumerate}

\item la fibre générique $(\!Z_\!U)_K$ de $\!Z_\!U$ est de codimension au moins 1 dans $X_K$.

\item Le morphisme de $k$-schémas formels $\G(\!Z_\!U)\hookrightarrow\G(\!U)$ est fortement bijectif (\textit{cf} \S\ref{fi}).

\end{enumerate}

\end{prop}

\begin{proof}
La question étant locale en $\!X$, on peut supposer que $\!X:=\Spf A$ est affine défini dans $R\{X_1,\ldots,X_N\}$ par un idéal $\!I$. Soit $I_K:=\!I\times_R K$. On peut supposer que $\!X$ est irréductible. En particulier, ceci entra\^\i ne que l'idéal $I_K$ est premier.

\textit{Premier cas: supposons que $A_K:=A\otimes_R K$ est analytiquement séparable}. Dans l'article \cite{sch}, Schappacher a montré que le lieu singulier de $X_K$ est contenu srtictement dans $X_K$. En particulier, ceci implique que la dimension de $(\!X_{\mathrm{sing}})_K$ est strictement inférieure à $d$ (car $(\!X_{\mathrm{sing}})_K$ est un fermé de $X_K$ et $X_K$ est irréductible.

\textit{Deuxième cas: supposons que $A_K:=A\otimes_R K$ n'est pas analytiquement séparable}. Par un argument d'analyse ultramétrique utilisé par Schappacher dans \cite{sch}, l'hypothèse entra\^\i ne qu'il existe $g\in K\{X_1,\ldots,X_N\}$ vérifiant les propriétés suivantes: $g\not\in I_K$ et $g(\varphi)=0$ pour tout $\varphi\in\G(\!X)(k)$. Quitte à multiplier par une puissance de $\pi$, on peut supposer en outre que $g\in R\{X_1,\ldots,X_N\}$. Le sous-$R$-schéma formel fermé défini par l'idéal $\!I + g$ est une solution au problème.

\end{proof}

\begin{lem}
\label{morphisme etale}
Soit $h:\!Y\rightarrow\!X$ un $\mathbb{D}$-morphisme étale de schémas formels sttf et réduits. Alors le $k$-schéma $\G_n(\!Y)$ (resp. $\G(\!Y)$) est (canoniquement) isomorphe au $k$-schéma $\G_n(\!X)\times_{X_0} Y_0$ (resp. $\G(\!X)\times_{X_0} Y_0)$), pour tout $n\in\mathbf{N}$.
\end{lem}

\begin{proof}
Dans le cas où $\!X$ est le complété $\pi$-adique d'une variété sur $k$, cette preuve est également faite dans \cite{kv}. On peut supposer $\!X$ et $\!Y$ affines. Soit $S=\Spec A$ un $k$-schéma affine. La donnée d'un $k$-morphisme $S\rightarrow\G_n(\!X)\times_{X_0} Y_0$ équivaut à celle d'un diagramme commutatif de $k$-morphismes
$$
\xymatrix{S\ar[r]\ar[d]&\G_n(\!X)\ar[d]^{\theta_0^n}\\
Y_0\ar[r]^{h_0}&X_0}
$$
Par adjonction, celle-ci équivaut à la donnée du diagramme commutatif de $R_n$-morphismes
$$
\xymatrix{\Spec(L(A))\otimes_{L(k)}R_n)\ar[rr]& & X_n &\\
&Y_n\ar@{-->}[ur]^{h_n}&&X_0\ar[ul]\\
S=\Spec(L(A(S))\otimes_{L(k)}R_0)\ar[uu]^{u}\ar[rr]\ar@{-->}[ur]&  & Y_0\ar@{-->}[ul]\ar[ur]^{h_0}}
$$
Comme le $R_n$-morphisme $u$ est une immersion nilpotente et que le $R_n$-morphisme $h_n$ est étale (\textit{cf} \S\ref{lissite}), il existe un unique $R_n$-morphisme $v:\Spec(L(A(S))\otimes_{L(k)}R_n)\rightarrow Y_n$, qui complète le diagramme ci-dessus en un diagramme de $R_n$-morphismes, où tous les triangles sont commutatifs
$$
\xymatrix{\Spec(L(A(S))\otimes_{L(k)}R_n)\ar@{-->}[dr]_{v}\ar[rr]& & X_n &\\
&Y_n\ar[ur]^{h_n}&&X_0\ar[ul]\\
S\ar[uu]^{u}\ar[rr]\ar[ur]&  & Y_0\ar[ul]\ar[ur]^{h_0}}
$$
Par adjonction, on déduit l'existence d'un $k$-morphisme de schéma $v:S\rightarrow\G_n(\!Y)$. L'unicité de ce morphisme découle de la construction.
\end{proof}


\subsection{Le cas lisse}


\begin{exe}
Si $\!X=\mathbb{B}_R^d$, le $k$-schéma $\G_n(\!X)$ est $k$-isomorphe à $\mathbf{A}_k^{(n+1)d}$.
\end{exe}

\begin{lem}
\label{le cas lisse}
Soit $\!X$ un $R$-schéma formel sttf lisse sur $R$ de dimension $d$. Alors on a les propriétés suivantes:
\begin{enumerate}

\item pour tout $n\in\mathbf{N}$ et tout $m\geq n$, les $k$-morphismes $\pi_{n,\!X}$ et $\theta^m_n$ sont surjectifs.

\item Pour tout $n\in\mathbf{N}$, le $k$-morphisme $\theta_n^{n+1}$ est une fibration localement triviale pour la topologie de Zariski de fibre $\mathbf{A}^d_k$.

\end{enumerate}
\end{lem}

\begin{proof}
L'assertion \ref{le cas lisse}/1 est une conséquence de \ref{le cas lisse}/2. Comme $\!X$ est lisse, pour tout $x\in X_0$, il existe un sous-$R$-schéma formel affine $\!U\hookrightarrow\!X$, ouvert contenant $x$, et un $R$-morphisme de schémas formels étale $\!U\rightarrow\mathbb{B}_R^d$. Le lemme \ref{morphisme etale} et le fait que $\G_n(\mathbb{B}_R^d)=\mathbf{A}_k^{(n+1)d}$ assurent le résultat.
\end{proof}


\subsection{L'interprétation des points géométriques}
\label{fi}


Soit $x$ un point géométrique de $\G(\!X)$ et $\overline{x}$ un point géométrique au-dessus de $x$, \textit{i.e.} un $k$-morphisme de schémas $\Spec \kappa(x)^{\mathrm{alg}}\rightarrow \G(\!X)$ avec $\kappa(x)^{\mathrm{alg}}$ une clôture algébrique du corps résiduel $\kappa(x)\supset k$. Par adjonction, il lui correspond un unique $R$-morphisme de schémas formels $\underline{x}:\Spf R_{\kappa(x)^{\mathrm{alg}}}\rightarrow \!X$. Posons $k':=\kappa(x)^{\mathrm{alg}}$. Désignons par $K'$ le corps des fractions de l'anneau $R_{k'}$. Par extension des scalaires, les ensembles
$$
\Hom_\mathbb{D}(\Spf R_{k'},\!X)\rightarrow\Hom_{\Spf R_{k'}}(\Spf R_{k'},\!X\hat{\times}_\mathbb{D}\Spf R_{k'})
$$
s'identifient canoniquement. On peut donc associer à $\overline{x}$ un (unique) rig-point d'une extension de $\!X$. En particulier, si $\!X$ est affine, un point géométrique de $\G(\!X)$ s'identifie à un (unique) point de l'espace sous-jacent d'une extension de $X_K$ par un corps $K'$ (\textit{cf} lemme 3.4 de \cite{bl1}).

Soit $h:\!Y\rightarrow\!X$ un $\mathbb{D}$-morphisme de schémas formels \textit{sttf}. On dira que $h:\G(\!Y)\rightarrow\G(\!X)$ est \textit{injectif} si l'application d'ensembles sous-jacente est injective. On dit que $h$ est \textit{fortement injectif} si, pour tout corps $F$ parfait (et séparable) contenant $k$, l'application canonique
$$
\G(\!Y)(F)\rightarrow\G(\!X)(F)
$$
est injective. Si $h$ est fortement injectif, alors $h$ est injectif. Soient $x$ et $y$ deux points \textit{distincts} de $\G(\!Y)$, de corps résiduels $\kappa(x)$ et $\kappa(y)$ dans $\G(\!Y)$. Soit $F$ le corps parfait et séparable défini de la manière suivante
$$
F:=(\mathrm{Frac}(\kappa(x)^{\mathrm{alg}}\otimes_{k^{\mathrm{alg}}}\kappa(y)^{\mathrm{alg}}))^{\mathrm{alg}}
$$
où $\mathrm{Frac}(A)$ est le corps des fractions de $A$, si $A$ est un anneau intègre. Le corps $F$ contient les corps $\kappa(x)$ et $\kappa(y)$ (donc $k$). Soit $\varphi_x:\Spec (\kappa(x)^{\mathrm{alg}})\rightarrow\G(\!Y)$ (\textit{resp.} $\varphi_y:\Spec (\kappa(y)^{\mathrm{alg}})\rightarrow\G(\!Y)$) le $k$-morphisme correspondant à $x$ (\textit{resp.} à $y$). On a les diagrammes commutatifs suivants
$$
\xymatrix{&&\Spec F\ar[drr]^{p_2}\ar[dll]_{p_1}&&\\
\Spec (\kappa(x)^{\mathrm{alg}})\ar[rr]^{\varphi_x}\ar[drr]^{\varphi_1}&&\G(\!Y)\ar[d]^h&&\Spec (\kappa(y)^{\mathrm{alg}})\ar[ll]_{\varphi_y}\ar[dll]_{\varphi_2}\\
&&\G(\!X)&&}
$$
Par hypothèse et construction, les morphismes $\varphi_x\circ p_1$ et $\varphi_x\circ p_2$ sont distincts. L'hypothèse d'injectivité entraine alors que les morphismes $\varphi_1\circ p_1$ et $\varphi_2\circ p_2$ sont distincts. En particulier, leurs images, \textit{i.e.} les points $h(x)$ et $h(y)$ de $\G(\!X)$, sont également distincts.

De même, on dit que $h$ est \textit{surjectif} (\textit{bijectif}) si l'application d'ensembles sous-jacente est surjective (bijective). On dit que $h$ est \textit{fortement surjectif} (\textit{fortement bijectif}) si, pour tout corps $F$ parfait (et séparable) contenant $k$, l'application canonique
$$
\G(\!Y)(F)\rightarrow\G(\!X)(F)
$$
est surjective (bijective). Si $h$ est fortement surjectif (fortement bijectif), alors $h$ est surjectif (bijectif). En effet, soit $x\in\G(\!X)$. Soit $\tilde{x}\in\G(\!X)(F)$ un $F$-point géométrique de $\G(\!X)$ au-dessus de  $x$. Alors, par hypothèse, il existe un $F$-point $\tilde{y}$ de $\G(\!Y)$ qui est un antécédant de $\tilde{x}$ par $h$. En particulier, le point $\tilde{y}(F)$ de $\G(\!Y)$ est un antécédent de $x$ par $h$.

\begin{exe}
Soit $\!X$ un $R$-schéma formel \textit{sttf}, génériquement lisse. Il possède donc un modèle de Néron faible $\sigma:\!U\rightarrow\!X$ sttf (\textit{cf} \cite{bs}). Par définition, le $R$-morphisme $\sigma$ est fortement injectif et fortement surjectif.
\end{exe}

Pour toutes ces raisons, nous identifierons, dans la suite de l'article, toutes les fois où cela sera nécessaire, les notions de points de $\G(\!X)$ et de $F$-points de $\G(\!X)$, pour $F$ une extension parfaite de $k$. Ainsi, si $x\in\G(\!X)$, on notera $\varphi_x:\Spf R_{k'}\rightarrow\!X$ (ou plus simplement $\varphi$) le $\mathbb{D}$-morphisme correspondant à un point géométrique de $\G(\!X)$ au-dessus de $x$.


\section{Cylindres et leur mesure}



\subsection{Les anneaux de Grothendieck.}


Soit $k$ un corps. On note $\!M:=K_0(\underline{Var}_{/k})$ le groupe abélien engendré par les symboles $[S]$, pour $S$ une variété sur $k$ (\textit{i.e.} un $k$-schéma de type fini réduit et séparé), avec les relations $[S]=[S']$ si $S$ et $S'$ sont isomorphes et $[S]=[S'] +[S\backslash S']$ si $S'$ est une sous-variété fermée de $S$. Il existe une structure naturelle d'anneau sur $\!M$, le produit étant induit par le produit fibré. La classe du point est l'élément neutre de cet anneau. On la notera $\mathbf{1}$. 

À tout ensemble constructible $C$ d'une variété $S$, on peut associer naturellement une classe $[C]$ dans $\!M$. Si $C$ et $C'$ sont deux ensembles constructibles d'une variété $S$, on a la relation $[C\cup C']=[C]+[C']-[C\cap C']$. 

On notera $\mathbf{L}$ la classe de $\mathbf{A}_k^1$ dans $\!M$, \textit{i.e.} $\mathbf{L}:=[\mathbf{A}_k^1]$. On désignera par $\!M_{\mathrm{loc}}:=\!M[\mathbf{L^1}]$ le localisé de $\!M$ par rapport au système multiplicatif $\{\mathbf{1},\mathbf{L}, \mathbf{L}^2,\ldots\}$.

Soit $F^m\!M_{\mathrm{loc}}$ le sous-groupe de $\!M_{\mathrm{loc}}$ engendré par les $[S]\mathbf{L}^{-i}$ tels que $\dim S-i\leq -m$, et $\widehat{\!M}$ le séparé, complété de $\!M_{\mathrm{loc}}$ suivant cette filtration. On notera $F^\bullet$ cette filtration. On note $\overline{\!M_{\mathrm{loc}}}$ l'image de $\!M_{\mathrm{loc}}$ dans $\widehat{\!M}$.

La filtration $F^\bullet$ définit une topologie métrisable sur $\widehat{\!M}$. On note
$$
\xymatrix{\parallel\hskip 1mm\parallel:\widehat{\!M}\ar[r]&\mathbf{R}_{\geq 0}}
$$
l'application  telle que
$$
\lVert{a}\rVert=
\begin{cases}
2^{-n}& \text{si $a\in F^n\widehat{\!M}$ et $a\not\in F^{n+1}\widehat{\!M}$},\\
0 & \text{$a=0$}. 
\end{cases}
$$
Cette application est la norme induite par cette filtration. Elle munit $\widehat{\!M}$ d'une structure d'anneau normé non archimédien.

Soit $S$ une $k$-variété. Tout sous-ensemble constructible $C$ de $S$ peut s'écrire comme réunion finie disjointe de sous-$k$-variétés de $S$, \textit{i.e.} il existe des sous-$k$-variétés $(S_i)_{1\leq i\leq n}$ de $S$ telles que
$$
C=\sqcup_{1\leq i\leq n} S_i.
$$
On définit alors la \textit{dimension} de $C$ comme la borne supérieure des dimensions des $S_i$ pour $1\leq i\leq n$. Il est clair que cette définition est indépendante de la partition choisie.


\subsection{Les morphismes par morceaux}


\begin{defi}
Soient $X$ et $Y$ deux $k$-sch\'emas, $A$ (resp. $B$) une partie constructible de $X$ (resp. $Y$). On dit qu'une application $\pi:A\rightarrow B$ est un morphisme par morceaux s'il existe une partition finie de $A$ en sous-schémas localement fermés de $X$ telle que la restriction de $\pi$ à chacun de ces sous-schémas soit induite par un morphisme de $k$-schémas.
\end{defi}

\begin{defi}
Avec les notations ci-dessus, on dit que l'application $\pi:A\rightarrow B$ est une fibration triviale par morceaux de fibre $F$, s'il existe une partition finie de $B$ par des sous-ensembles $S$ localement ferm\'es dans $Y$ tel que $\pi^{-1}(S)$ est localement ferm\'e dans $X$ et isomorphe, en tant que $k$-schéma, \`a $S\times_k F$, $\pi$ correspondant, par cet isomorphisme \`a la projection $S\times_k F\rightarrow S$.
\end{defi}

\begin{defi}
Avec les m\^emes notations, on dit que l'application $\pi$ est une fibration triviale par morceaux sur $C$, avec $C$ un sous-ensemble constructible de $B$, si la restriction de $\pi$ \`a $\pi^{-1}(C)$, \textit{i.e.} $\pi_{\mid \pi^{-1}(C)}:\pi^{-1}(C)\rightarrow C$, est une fibration triviale par morceaux.
\end{defi}

\begin{rem}
Si $\pi:A\rightarrow B$ est une fibration triviale par morceaux de fibre $F$, les lois algébriques dans $\!M$ assurent que $[A]=[B].[F]$.
\end{rem}


\subsection{Les cylindres.}


Fixons un $R$-schéma formel \textit{sttf}.

\subsubsection{La définition et les premiers exemples.}

\begin{defi}
On dit qu'une partie $A$ de $\G(\!X)$ est un sous-ensemble cylindrique de rang $n$ de $\G(\!X)$ (ou plus simplement un cylindre de rang $n$ de $\G(\!X)$), si $A=\pi_{n,\!X}^{-1}(C)$, o\`u $C$ d\'esigne une partie constructible de $\G_n(\!X)$. On dit que $A$ est un cylindre si $A$ est un cylindre  d'un certain rang $n$.
\end{defi}

\begin{defi}\label{pro-cylindre}
On dit qu'une partie $A$ de $\G(\!X)$ est un pro-cylindre (resp. un ind-cylindre) si elle est intersection (resp. r\'eunion) d\'enombrable de cylindres. 
\end{defi}

\begin{lem}
Soit $A$ un cylindre de rang $n$ de $\G(\!X)$. Alors $A=\pi_{n,\!X}^{-1}(\pi_{n,\!X}(A))$.
\end{lem}

\begin{proof}
Évident, car $A:=\pi_{n,\!X}^{-1}(A_n)$ avec $A_n$ ensemble constructible de $\G_n(\!X)$.
\end{proof}

\begin{exes}

\begin{enumerate}

\item L'espace $\G(\!X)$ est un 0-cylindre, puisque $\pi_{0,\!X}^{-1}({(X_0)}_{\mathrm{red}})=\G(\!X)$.

\item Soit $\!U\hookrightarrow\!X$ un sous-$R$-sch\'ema formel ouvert de $\!X$ quasi-compact. Le $k$-sch\'ema $\G(\!U)\hookrightarrow\G(\!X)$ est un 0-cylindre de $\G(\!X)$ de rang 0. Plus pr\'ecis\'ement, on a l'\'egalit\'e:$$\G(\!U)=\pi_{0,\!X}^{-1}(\G_0(\!U)).$$ En effet, par d\'efinition on a que $\pi_{0,\!X}(\G(\!U))\subset \G_0(\!U)$. Soit $\varphi\in\pi_{0,\!X}^{-1}(\G_0(\!U))$. Ceci entra\^\i ne en particulier que l'image de $\varphi:\Spf R_{k'}\rightarrow \!X$ est contenue dans $\!U$. Soit $\!V$ un ouvert affine de $\!U$ (donc ouvert de $\!X$) contenant l'image topologique de $\varphi$. Le $R$-morphisme $\varphi$ se factorise par $\!V$ donc par $\!U$.

\item Soit $\!Z\hookrightarrow\!X$ un sous-$R$-sch\'ema ferm\'e de $\!X$, muni du morphisme canonique. Alors $\G(\!Z)$ est un pro-cylindre. En effet, on a l'\'egalit\'e:$$\G(\!Z)=\cap_{n\in\mathbf{N}}\pi_{n,\!X}^{-1}(\G_n(\!Z)).$$
\end{enumerate}

\end{exes}

\subsubsection{Les propriétés essentielles.}

\begin{prop}
Tout cylindre est un sous-ensemble constructible de $\G(\!X)$. 
\end{prop}

\begin{proof}
Ce lemme découle du fait que l'image inverse d'un constructible est constructible (\textit{cf} \cite{ega} \S 7).
\end{proof}

\begin{lem}
\label{tout cylindre de rang n est de rang m}
Si $A$ est un cylindre de rang $n$, alors $A$ est aussi un cylindre de rang $m$, pour tout $m\geq n$.
\end{lem}

\begin{proof}
Ce lemme découle de la commutativité du diagramme suivant
$$
\xymatrix{\G(\!X)\ar[r]^{\pi_{m,\!X}}\ar[dr]_{\pi_{n,\!X}}&\G_m(\!X)\ar[d]^{\theta_n^m}\\
&\G_n(\!X)}
$$
et du fait que l'image inverse d'un constructible est constructible (\textit{cf} \cite{ega} \S 7).
\end{proof}

On note $\mathbf{C}_\!X$ l'ensemble des cylindres de $\G(\!X)$.

\begin{prop}
\label{cylindres sont un anneau booleen}
L'ensemble $\mathbf{C}_\!X$ des cylindres de $\G(\!X)$ est un anneau booléen. Autrement dit,

\begin{enumerate}

\item les ensembles $\G(\!X)\in\mathbf{C}_\!X$ et $\emptyset\in\mathbf{C}_\!X$.

\item L'ensemble $\mathbf{C_\!X}$ est stable par intersection finie.

\item L'ensemble $\mathbf{C}_\!X$ est stable par r\'eunion finie.

\item L'ensemble $\mathbf{C_\!X}$ est stable par passage au compl\'ementaire.
\end{enumerate}

\end{prop}

\begin{proof}
Grâce au lemme \ref{tout cylindre de rang n est de rang m}, on se ramène au cas de cylindres de même rang. Cette proposition découle alors du fait que les parties constructibles d'un ensemble vérifient ces propriétés et de la stabilité de ces propriétés pour l'image inverse (\textit{cf} \cite{ega} \S 7).
\end{proof}

\begin{lem}
\label{quasi-compacite de la topologie constructible}
Soit $(A_i)_{i\in I}$ une famille d\'enombrable de cylindres de $\G(\!X)$. Si $A:=\cup\displaystyle_{i\in I} A_i$ est \'egalement un cylindre, alors il existe un sous-ensemble $J$ fini de $I$ tel que $A:=\cup\displaystyle_{i\in J} A_i$.
\end{lem}

\begin{proof}
En remarquant que $\G(\!X)$ est quasi-compact (obtenu par recollement d'un nombre fini d'affines), ce lemme est une conséquence de la quasi-compacité de la topologie constructible (\textit{cf} \cite{ega} \S 7).
\end{proof}

Le résultat suivant, qui est une formulation du théorème de Schappacher \cite{sch}, est l'analogue du théorème de Greenberg \cite{g2} dans le cadre de la théorie des schémas formels.

\begin{thm}
\label{theoreme de Greenberg}
Soit $R$ un anneau de valuation discrète complet et $\!X$ un $R$-schéma formel sttf. Pour tout $n\geq 0$, il existe un entier $m_\!X(n)\geq n$, ne dépendant que de $\!X$, tel que, pour tout corps parfait $F$ contenant $k$, et tout $x\in \!X(R_F/(\pi)^{m_\!X(n)})$, l'image de $x$ dans $\!X(R_F/(\pi)^{n})$ peut être relevée en $\tilde{x}\in\!X(R_F)$. On appelle fonction de Greenberg du schéma formel $\!X$ l'application qui à $n\in\mathbf{N}$ associe le plus petit des entiers $m_\!X(n)$ définis ci-dessus. On notera $\gamma_\!X$ cette application.
\end{thm}

Ce théorème permet de démontrer le lemme suivant, qui joue un rôle crucial dans la définition de la mesure motivique pour les cylindres.

\begin{lem}
\label{1}
L'image $\pi_{n,\!X}(\G(\!X))$ de $\G(\!X)$ dans $\G_n(\!X)$ est un ensemble constructible de $\G_n(X_n)$. Plus généralement, si $A$ est un cylindre de rang $m$ de $G(\!X)$, alors $\pi_{n,\!X}(A)$ est un ensemble constructible de $\G_n(\!X)$, pour tout $n\geq 0$.
\end{lem}

\begin{proof}
Il nous suffit de démontrer la seconde assertion. Par définition, il existe un ensemble constructible $C_m$ de $\G_m(\!X)$ tel que $A=\pi_m^{-1}(C_m)$ (on omettra l'indice $\!X$ dans $\pi_{n,\!X}$). On peut supposer que $m=n$. Si $n\geq m$, $A=\pi_n^{-1}(({\theta_m^n})^{-1}(C_m))$. Si $m\geq n$, $\pi_n(A)=\theta_n^m(\pi_m(A_m))$; donc, par un théorème de Chevalley (\textit{cf} \cite{ega} \S 7), si $\pi_m(A)$ est constructible, $\pi_n(A)$ l'est aussi.

Supposons que $m=n$. Le théorème de Greenberg \ref{theoreme de Greenberg} assure que
$$
\pi_n(A)=\theta_n^{\gamma_\!X(n)}(({\theta_n^{\gamma_\!X(n)}})^{-1}(C_n)).
$$
Le théorème de Chevalley déjà cité (\textit{cf} \cite{ega} \S 7) permet alors de conclure.
\end{proof}

Les lemmes précédents permettent de borner la dimension des fibres des morphismes $\theta_n^m$.

\begin{lem}
\label{dimension de la fibre}
Soit $\!X$ un sch\'ema formel sttf sur $R$ dont la fibre générique $X_K$ est de dimension $d$. Alors:

\begin{enumerate}

\item pour tout $n\in\mathbf{N}$, $\dim\pi_{n,\!X}(\G(\!X))\leq(n+1)d$.

\item Pour tout $n,m\in\mathbf{N}$, tel que $m\geq n$, les fibres de $\pi_{m,\!X}(\G(\!X))\rightarrow\pi_{n,\!X}(\G(\!X))$ sont de dimension inf\'erieure \`a $(m-n)d$.
\end{enumerate}

\end{lem}

\begin{proof}
Supposons connue l'assertion \ref{dimension de la fibre}/2. Appliquons \ref{dimension de la fibre}/2 aux entiers $n\geq 0$ et $0$. On a $\pi_0(\G(\!X))\subset X_0$ et donc
$$
\dim \pi_n(\G(\!X))\leq nd +\dim \pi_0(\G(\!X))\leq (n+1)d.
$$
Pour montrer l'assertion \ref{dimension de la fibre}/2, il nous suffit de prouver que, pour tout $n\in\mathbf{N}$, la dimension de chaque fibre du morphisme $\pi_{n+1}(\G(\!X))\rightarrow\pi_n(\G(\!X))$ est inférieure ou égale à $d$. 

On peut supposer que $\!X$ est affine de la forme $\Spf R\{x_1,\ldots,x_N\}/(f_1,\ldots,f_m)$. Posons $f:=(f_1,\ldots,f_m)$. Soit $x\in\G(\!X)$. Soit $\tilde{\xi}\in\G(\!X)(\kappa(x)^{\mathrm{alg}})$ le $k$-morphisme qui correspond à un point géométrique au-dessus de $x$ et $\xi\in R^N$ le $N$-uplet de $R^N$ qui correspond à ce morphisme. On peut supposer que $k=\kappa(x)^{\mathrm{alg}}$ et que $\varphi\in R^N$. Soit $\!Y$ le $R$-schéma formel affine \textit{ttf} défini par le système de séries formelles restreintes $g(y)=f(\xi+\pi^{n+1}y)$, avec $y:=(y_1,\ldots,y_N)$ dans $\Spf R\{y\}$. Le $K$-morphisme d'algèbres topologiques qui à $x_i$ associe  $\xi_i + \pi^{n+1}y_i$ pour tout $1\leq i\leq N$ induit un isomorphisme d'espaces rigides entre $\!X_{\mathrm{rig}}$ et $\!Y_{\mathrm{rig}}$. En particulier, $\!Y_{\mathrm{rig}}$ est de dimension $d$. Il existe un $R$-modèle $\!Y'$ de $\!Y_{\mathrm{rig}}$ défini, dans $\mathbb{B}^N_R$, par le système d'équations $g(y)/\pi^{n+1}$. Tout point de la fibre est alors contenu dans la fibre spéciale $Y'_0$ de $\!Y'$. Le théorème 1 établi par Oesterlé dans \cite{o1} (\textit{cf} théorème \ref{theoreme de Osterle}) assure alors que la fibre est de dimension au plus $d$.
\end{proof}

Signalons cette conséquence du lemme \ref{dimension de la fibre}.

\begin{lem}
\label{nullite du lieu singulier}
\label{3}
Soit $\!X$ un $R$-sch\'ema formel sttf dont la fibre générique $X_K$ est de dimension $d$ et soit $\!S$ un sous-$R$-sch\'ema formel fermé, dont la fibre générique $S_K$ est de dimension strictement inf\'erieure \`a $d$. Soit $\gamma_\!S$ la fonction de Greenberg pour $\!S$. Alors pour tous $n,i,e\in\mathbf{N}$, tels que $n\geq i\geq \gamma_\!S(e)$, $\pi_{n,\!X}(\pi_{i,\!X}^{-1}\G_e(\!S))$ est de dimension inf\'erieure \`a $(n+1)d-e-1$. Autrement dit, pour tous $n\geq i\geq \gamma_\!S(e)$,
$$
[\pi_{n,\!X}(\pi_{i,\!X}^{-1}(\G_i(\!S)))].\mathbf{L}^{-(n+1)d}\in F^{e+1}\!M_{\mathrm{loc}}.
$$
\end{lem}

\begin{proof}
On oubliera dans cette preuve les indices $\!X$ et $\!S$. On peut supposer que $i=\gamma(e)$. Par le lemme \ref{dimension de la fibre}/2 appliqué à la projection:
$$
\pi_n(\pi_{\gamma(e)}^{-1}(\G_{\gamma(e)}(\!S)))\rightarrow\pi_e(\pi_{\gamma(e)}^{-1}(\G_{\gamma(e)}(\!S)))
$$
on obtient l'inégalité:
$$
\dim \pi_n(\pi_{\gamma(e)}^{-1}(\G_{\gamma(e)}(\!S)))\leq (n-e)d + \dim \pi_e(\pi_{\gamma(e)}^{-1}(\G_{\gamma(e)}(\!S))).
$$

Par ailleurs, par définition de la fonction de Greenberg, 
$$
\pi_e(\pi_{\gamma(e)}^{-1}(\G_{\gamma(e)}(\!S)))=\pi_e(\G(\!S)).
$$
Le lemme \ref{dimension de la fibre}/1, assure que $\dim \pi_e(\G(\!S)\leq (e+1)(d-1)$. Le résultat découle du fait que $(n-e)d+(e+1)(d-1)=(n+1)d-e-1$.

\end{proof}

\begin{cor}
Soit $\!X$ un $R$-schéma formel admissible. Soit $B$ un cylindre de $\G(\!X)$ de rang $m$ tel qu'il existe un sous-$R$-schéma formel $\!Z$ de $\!X$, dont le $k$-schéma de Greenberg associé contient $B$. Alors, pour tous $n\geq i\geq \max(m,\gamma_\!S(e))$,
$$
[\pi_{n,\!X}(B)].\mathbf{L}^{-(n+1)d}\in F^{e+1}\!M_{\mathrm{loc}}.
$$
\end{cor}

\begin{rem}
Soit $\!X$ un $R$-schéma formel admissible et réduit. Alors la proposition \ref{schappa}, le lemme \ref{3} et son corollaire vont entra\^\i ner que soit le volume du lieu singulier de $\!X$ est nul, soit le volume de $\G(\!X)$ est nul. Ce dernier cas étant sans grand intérêt, nous ne le traiterons pas explicitement dans les définitions et les preuves qui vont suivre, bien que celles-ci soient encore valides dans ce cas.
\end{rem}

\begin{exe}
Soit $e\in\mathbf{N}$. Un exemple important de cylindre est le sous-ensemble $\G^{(e)}(\!X)$ de $\G(\!X)$ défini de la manière suivante: on pose
$$
\G^{(e)}(\!X):=\G(\!X)\backslash(\pi_{e,\!X}^{-1}(\G_e(\!X_{\mathrm{sing}}))).
$$
Cette définition induit la décomposition ensembliste:
$$
\G(\!X)=(\bigcup_{e\in\mathbf{N}}\G^{(e)}(\!X))\sqcup \G(\!X_{\mathrm{sing}}).
$$
Soit $x\in\G(\!X)$. On pose $\varphi_x:\Spec \kappa(x)^p\rightarrow \G(\!X)$ le $k$-morphisme correspondant, avec $\kappa(x)^p$ une clôture parfaite du corps résiduel $\kappa(x)$ de $x$ dans $\G(\!X)$. Si $x\in\G^{(e)}(\!X)$, alors le $K$-morphisme $(\varphi_x)_{\mathrm{rig}}$ se factorise par le lieu lisse de la fibre générique $\!X_{\mathrm{rig}}$ de $\!X$ (cf \cite{bl3} \S 2 et 3) et le défaut de lissité de Néron de l'adjoint du morphisme $\varphi_x$ est inférieur ou égal à $e$ (cf \cite{bs} \S 3 et \cite{blr} \S 3.3 et lemme 3.3/2).

\end{exe}

\subsubsection{La notion de cylindre stable.}

\begin{defi}
On dit qu'un sous-ensemble $A$ de $\G(\!X)$ est un cylindre stable de rang $n$ si:

\begin{enumerate}

\item l'ensemble $A$ est un cylindre de rang $n$.

\item Pour tout $m\geq n$, la restriction du morphisme de transition $\pi_{m+1}(\G(\!X))\rightarrow\pi_{m}(\G(\!X))$ est une fibration triviale par morceaux de fibre $\mathbf{A}_k^d$ sur $\pi_m(A)$.
\end{enumerate}
\end{defi}

Soit $\mathbf{C}_{0,\!X}$ la famille des sous-ensembles cylindriques de $\G(\!X)$ qui sont stables pour un certain rang $n$.

\begin{prop}
L'ensemble $\mathbf{C}_{0,\!X}$ des cylindres stables de $\G(\!X)$ est un idéal de $\mathbf{C_\!X}$. Autrement dit, il vérifie les propriétés suivantes:

\begin{enumerate}

\item l'ensemble $\mathbf{C}_{0,\!X}\subset \mathbf{C}_{\!X}$ et contient $\emptyset$.

\item Si $A,B\in\mathbf{C}_{0,\!X}$ sont deux cylindres stables disjoints, alors leur r\'eunion est un cylindre stable.

\item Si $A,B\in\mathbf{C}_{0,\!X}$ sont deux cylindres stables, alors leur r\'eunion est un cylindre stable.

\item Si $A\in\mathbf{C}_{0,\!X}$ et $B\in\mathbf{C}$, alors $A\cap B$ est un cylindre stable.
\end{enumerate}

\end{prop}

\begin{proof}
Ces propriétés découlent de la proposition \ref{cylindres sont un anneau booleen} et de la définition des cylindres stables.
\end{proof}

\begin{exe}
Le cylindre $\G^{(e)}(\!X)$ est stable (\textit{cf} lemme \ref{existence de la mesure}).
\end{exe}

L'intérêt des cylindres stables est que l'on peut définir une mesure na\"{i}ve sur cet ensemble de la manière suivante.

\begin{prop}
\begin{enumerate}
\item Il existe un unique morphisme additif
$$
\mu_{0,\!X}:\mathbf{C}_{0,\!X}\longrightarrow \!M
$$
tel que $\mu_{0,\!X}(A)=[\pi_{n,\!X}(A)]\mathbf{L}^{-(n+1)d}$, pour tout cylindre $A$ stable au rang $n$.

\item L'application $A\mapsto{\mu}_{0,\!X}(A)$ est $\sigma$-additive sur $A\in\mathbf{C}_{0,\!X}$.

\item Pour $A$ et $B$ dans $A\in\mathbf{C}_{0,\!X}$, $\lVert{\mu}_{0,\!X}(A\cup B)\lVert\leq \max(\lVert{\mu}_{0,\!X}(A)\rVert,\lVert{\mu}_{0,\!X}(B)\rVert)$. Si $A\subset B$, $\lVert{\mu}_{0,\!X}(A)\rVert\leq\lVert{\mu}_{0,\!X}(B)\rVert.$

\end{enumerate}
\end{prop}

\begin{proof}

Soit $A$ un cylindre stable de rang $n$ de $\G(\!X)$. Comme, par définition, l'application $\theta_n^m$ est sur $\pi_n(A)$ une fibration localement triviale par morceaux de fibre $\mathbf{A}^d_k$ pour tout $m\geq n$, on a dans ${\!M}_{\mathrm{loc}}$ la relation:
$$
[\pi_{m,\!X}(A)]=[\pi_{n,\!X}(A)]\mathbf{L}^{(m-n)d}.
$$
Cette relation permet de donner un sens à la définition, puisque $A$ est un cylindre stable également au rang $m\geq n$.

Soit $A:=\cup_{i\in\mathbf{N}} A_i$ un cylindre de $\G(\!X)$ stable au rang $n$, réunion de cylindres stables disjoints deux à deux. L'additivité de l'application $[\hskip 1mm]$ entra\^\i ne l'additivité de $\mu_{0,\!X}$. Pour la $\sigma$-additivité, on se ramène au cas de l'additivité grâce au lemme \ref{quasi-compacite de la topologie constructible}.
\end{proof}

\subsubsection{La mesure pour les cylindres.}

\begin{prop}
\label{proprietes de la mesure}
Soit $\!X$ un $R$-schéma formel admissible et réduit.
\begin{enumerate}

\item Pour tout cylindre $A\in\mathbf{C}_\!X$, la limite
$$
\tilde{\mu}_\!X(A):=\lim_{e\rightarrow\infty}\mu_{0,\!X}(A\cap\G^{(e)}(\!X))
$$
existe dans $\widehat{\!M}$

\item Si $A\in\mathbf{C}_{0,\!X}$, alors $\tilde{\mu}_\!X(A)$ et $\mu_{0,\!X}(A)$ co\"{i}ncident dans $\!M_{\mathrm{loc}}$.

\item L'application $A\mapsto\tilde{\mu}_\!X(A)$ est $\sigma$-additive sur $A\in\mathbf{C}_\!X$.

\item Pour $A$ et $B$ dans $A\in\mathbf{C}_\!X$, $\lVert\tilde{\mu}_\!X(A\cup B)\rVert\leq \max(\lVert\tilde{\mu}_\!X(A)\rVert,\lVert\tilde{\mu}_\!X(B)\rVert)$. Si $A\subset B$, $\lVert\tilde{\mu}_\!X(A)\rVert\leq\lVert\tilde{\mu}_\!X(B)\rVert$.

\end{enumerate}
\end{prop}

La preuve s'appuie sur le lemme clé suivant.

\begin{lem}
\label{existence de la mesure}
\label{2}

Soit $\!X$ un $R$-schéma formel admissible et réduit de pure dimension $d$. Il existe $c_\!X\in\mathbf{N}\backslash\lbrace 0\rbrace$ tel que, pour tout $e$ et tout $n\in\mathbf{N}$ avec $n\geq c_\!Xe$, on ait le r\'esultat suivant:

\begin{enumerate}

\item l'application: 
$$
\theta_n:\pi_{n+1,\!X}(\G(\!X))\rightarrow\pi_{n,\!X}(\G(\!X))
$$ 
est une fibration triviale par morceaux sur $\pi_{n,\!X}(\G^{(e)}(\!X))$ de fibre $\mathbf{A}_k^d$.

\item En particulier,
$$
\lbrack\pi_{n,\!X}(\G_e(\!X))\rbrack=\lbrack\pi_{c_\!Xe,\!X}(\G_e(\!X))\rbrack\mathbf{L}^{d(n-c_\!Xe)}.
$$

\end{enumerate}

\end{lem}

\begin{proof}[Démonstration du lemme \ref{existence de la mesure}]

Le point \ref{existence de la mesure}/2 est une conséquence directe de l'assertion \ref{existence de la mesure}/1. On peut supposer que $\!X$ est affine de la forme 
$$
\!X=\Spf R\{x_1,\ldots,x_N\}/\!I\hookrightarrow\mathbb{B}_R^N.
$$ 
Soit $\!J$ l'idéal de $R\{x_1,\ldots,x_N\}$ engendré par les produits $h\delta$, où $h$ appartient au conducteur de $\!I$ dans l'idéal de $R\{x_1,\ldots,x_N\}$ défini par $m\leq N$ générateurs $g_1,\ldots,g_m\in \!I$ et où $\delta$ est un mineur $m\times m$ de la matrice jacobienne $(\partial g_j/\partial x_i)_{1\leq j\leq m, 1\leq i\leq N}$. Le lieu singulier est alors défini par l'idéal $\sqrt{\!I+\!J}$ (\textit{ref} lemme \ref{autre definition du lieu singulier}). Soit $c_\!X$ le plus petit entier naturel tel que la puissance $c_\!X$-ème radical de l'idéal $\!I+\!J$ soit contenue dans l'idéal $\!I+\!J$.

\textit{On peut supposer que $\!X\hookrightarrow\mathbb{B}_R^N$ est d'intersection complète}. En effet, soit $x$ un point de $\G(\!X)$. Soit $\varphi:\Spf R_{k'}\rightarrow\!X$ le $\mathbb{D}$-morphisme qui correspond à un point géométrique au-dessus de $x$. Supposons que $x\not\in \pi_{\!X,e}^{-1}(\G_e(\!X))$. Alors il existe $g\in \sqrt{\!I+\!J}$ tel que
$$
g(\varphi)\not\equiv 0 \ \ \mod \pi^{e+1}. \leqno (0)
$$
En particulier, la congruence $(0)$ entra\^\i ne l'existence d'un élément $u$ de l'idéal $\!J$ tel que
$$
u(\varphi)\not\equiv 0 \ \ \mod \pi^{c_\!Xe+1}.
$$
On peut supposer $u=h\delta$, où $h\in((g_1,\ldots,g_m):\!I)$ et $\delta$ est un mineur de la matrice jacobienne $(\partial g_j/\partial x_i)_{1\leq j\leq m, 1\leq i\leq N}$. On peut supposer que $m=N-d$. En effet, $\varphi_e\in R^N_e$ appartient à $\G_e(\!X_{\mathrm{sing}})$ si et seulement si $q(\varphi_e)=0$ dans $R_e$ pour tout $q\in\sqrt{\!I+\!J}R_e$. Ceci reste vrai, en particulier, pour les éléments $q$ de la forme $h\delta$ ci-dessus avec $m=N-d$.

Le cylindre $\G^{(e)}(\!X)$ peut donc être recouvert par un nombre fini de cylindres de $\G(\mathbb{B}^N_R)$ de la forme:
$$
A:=\{\varphi\in\G(\mathbb{B}^N_R)\mid \ \ (h\delta)(\varphi)\not\equiv 0 \mod \pi^{ce+1}\}.
$$
Il nous suffit donc de prouver que, pour tout $n\geq c_\!Xe$, l'application $\theta_n$ est une fibration triviale par morceaux de fibre $\mathbf{A}_k^d$ sur $\pi_{n,\!X}(\G(\!X)\cap A)$ ($\G(\!X)\cap A$ est un cylindre de $\G(\!X)$). Comme $h(\varphi)\not=0$ pour tout $\varphi\in A$,
$$
\G(\!X)\cap A=\G(\Spf R\{x_1,\ldots,x_N\}/(f_1,\ldots,f_{N-d}))\cap A.
$$
On peut donc supposer que $\!X\hookrightarrow\mathbb{B}_R^d$ est d'intersection complète, \textit{i.e.} l'idéal $\!I$ est engendré par $N-d$ éléments $(f_1,\ldots,f_{N-d})$. Notons $\Delta$ la matrice jacobienne $(\partial f_j/\partial x_i)_{1\leq j\leq N-d, 1\leq i\leq N}$.

\textit{Montrons que la fibre de $\theta_n$ au-dessus d'un point quelconque $x$ de $\G(\!X)$ est $\kappa(x)$-isomorphe à $\mathbf{A}_{\kappa(x)}^d$}. Soit $x\in\G(\!X)$. Quitte à étendre les scalaires, on peut supposer que $k=\kappa(x)^{\mathrm{alg}}$. Soit alors $\varphi:\mathbb{D}\rightarrow \!X$ le $R$-morphisme correspondant au point $x$. Soit $e'\in\mathbf{N}$, $e'\leq c_\!Xe$. Posons
$$
A':=\{\varphi\in A\mid \ \ \ord_\pi\delta(\varphi)=e'\ \  \textrm{et}\ \  \ord_\pi\delta'(\varphi)\geq e'
$$
$$
\textrm{pour tout mineur $(N-d)\times(N-d)$ $\delta'$ de $\Delta$}\}.
$$
L'ensemble $A'$ est un cylindre de $\G(\!X)$. Il nous suffit donc de prouver que l'application $\theta_n$ est une fibration par morceaux triviale sur $\pi_{n,\!X}(\G(\!X)\cap A')$. On peut supposer, quitte à réordonner les $x_i$, que le mineur $\delta$ se calcule à partir des $N-d$ premières colonnes de $\Delta$. L'ensemble $\theta_n^{-1}(x)$ est décrit par les équations en $y\in R^N$
$$
f_j(\varphi+\pi^{n+1}y)=\pi^{n+1}D_\varphi f_j.y+\pi^{2n+2}F_j(y)=0
$$
pour tout $1\geq j\geq N-d$, où $F_j\in R\{x_1,\ldots,x_N\}$. Soit $M$ la matrice adjointe de la sous-matrice de $\Delta$ définie par les $N-d$ premières colonnes. En particulier, cette matrice à coefficients dans $R\{x_1,\ldots,x_N\}$ vérifie
$$
M.\Delta=(\delta I_{N-d}, W)
$$
où $W$ est une $d\times N-d$ matrice satisfaisant à
$$
W(\varphi)\equiv 0\mod(\pi^{e'}).
$$
En effet, soit $\widetilde{\Delta(\varphi)}$ la sous-matrice $N-d\times N-d$ de $\Delta(\varphi)$ formée des $N-d$ premières colonnes, que l'on note $\Delta_i(\varphi)$ pour $1\leq i\leq N-d$. Soit $j$ un entier naturel tel que $1\leq j\leq d$. Soit $W_j(\varphi)$ la $j$-ème colonne de $W(\varphi)$ (\textit{i.e.} la $N-d+j$-ème colonne de $M(\varphi).\Delta(\varphi)$). Le système
$$
\widetilde{\Delta(\varphi)}.X=W_j(\varphi), \leqno (0')
$$
où $X=(X_i)_{1\leq i\leq N-d}\in R^{N-d}$ est un vecteur colonne, admet, pour tout $1\leq j\leq d$, une solution, qui se décrit explicitement gr\^ace à la règle de Cramer par
$$
X_i=\mathrm{det}(\Delta_1(\varphi),\ldots,\Delta_{i-1}(\varphi),W_j,\Delta_{i+1}(\varphi),\ldots,\Delta_{N-d}(\varphi))/\mathrm{det}\Delta(\varphi)
$$
(qui est un élément de $R$ puisque le mineur $N-d\times N-d$ de $\Delta(\varphi)$ formé des $N-d$ premières colonnes est de valuation minimale parmi les mineurs $N-d\times N-d$ de $\Delta(\varphi)$).

La fibre est alors décrite par le système suivant
$$
\pi^{-e'}M.\Delta(\varphi) +\pi^{n+1-e'}(\ldots)=0. \leqno (1)
$$
Le lemme de Hensel (\textit{cf} \cite{bourb} III \S 4-3) assure alors que la fibre de $\theta_n$ au-dessus de $x$ est simplement l'ensemble des $y_0\in k^N$ vérifiant les équations linéaires
$$
\pi^{-e'}M.\Delta(y_0)\equiv 0\mod(\pi).
$$
La fibre $\theta_n^{-1}(x)$ est donc naturellement munie d'une structure de $\kappa(x)$-espace vectoriel, de dimension $d$, exprimant les $N-d$ premières coordonnées en fonction de combinaisons linéaires des $d$-dernières, dont les coefficients sont paramétrés par les premières coordonnées de $\varphi$ (\textit{cf} \S \ref{description de R}).

\textit{Nous allons conclure cette preuve par un argument ``noethérien''}. Soit $S_i$ un sous-$k$-schéma localement fermé de $\G_n(\!X)$ contenu dans $\pi_n(\G(\!X)\cap A')$. Le diagramme commutatif suivant
$$
\xymatrix{\pi_{n+1}(\G(\!X)\cap A')\ar[r]^{\theta_n} & \pi_n(\G(\!X)\cap A')\\
T_i:=\theta_n^{-1}(S_i)\ar[r]\ar@{^{(}->}[u] & S_i\ar@{^{(}->}[u]}
$$
assure qu'il nous suffit de démontrer qu'un $k$-morphisme $p:T\rightarrow S$ de schémas localement de type fini, dont les fibres vérifient la condition ci-dessus, est une fibration par morceaux triviale, \textit{i.e.} il existe une partition finie $(S_i)_i$ de $S$ en sous-$k$-schémas localement fermés telle que, pour tout $i$, le schéma $p^{-1}(S_i)$ est $k$-isomorphe à $S_i\times_k \mathbf{A}^d_k$. 

Le lemme \ref{existence ouvert} ci-dessous assure que, sous ces hypothèses, il existe un ouvert dense $U_0\hookrightarrow S$ de $S:=Z_0$ tel que $p^{-1}(U_0)$ soit $k$-isomorphe à $\mathbf{A}^d_k\times_k U_0$. Soit $Z_1$ le complémentaire de $U_0$ dans $S:=Z_0$. Nous pouvons alors appliquer à nouveau ce raisonnement au sous-$k$-schéma (localement fermé) $Z_1$ de $\G_n(\!X)$. En réitérant le procédé, on obtient une suite décroissante de fermés de $S$, $S:=Z_0\supset Z_1\supset\ldots$, qui, par noethérianité, stationne au rang $r$. La partition cherchée est alors donnée par les ouverts $U_i$ pour $0\leq i\leq r$.
\end{proof}

\begin{lem}
\label{existence ouvert}
Soient $S$ et $T$ deux $k$-schémas de type fini et $p:T\rightarrow S$ un $k$-morphisme de schémas tel que, pour tout $x\in S$, le $\kappa(x)$-schéma $p^{-1}(x)$ est $\kappa(x)$-isomorphe à $\mathbf{A}^d_{\kappa(x)}$. Alors il existe un ouvert dense $U$ de $S$ tel que $p^{-1}(U)$ est $U$-isomorphe à $\mathbf{A}^d_U$.
\end{lem} 

\begin{proof}
La propriété étant locale sur la base, on peut supposer que $S:=\Spec A$, où $A$ est une $k$-algèbre de type fini. Soit $(S_i)_{i\in I}$, $I$ fini, l'ensemble des composantes irréductibles de $S$. Pour tout $i\in I$, soit $\eta_i$ le point générique de $S_i$. Soit $D(g_i)$ un ouvert affine de $S$ contenant $\eta_i$. L'ouvert de $S$, $D(g_i)\backslash (\cup_{j\not= i} S_j)$, contient $\eta_i$ et est contenu dans $S_i$. Quitte à restreindre encore $S$, on peut supposer en outre que $S$ est irréductible. Désignons par $\eta$ son point générique.

\textit{Commen\c{c}ons par prolonger le $\kappa(\eta)$-morphisme de schémas $p^{-1}(\eta)\rightarrow\mathbf{A}_{\kappa(\eta)}^d$ à un morphisme au-dessus de $U$, pour $U$ un ouvert affine de $S$ de la forme $D(g)$ avec $g\in A$}. La donnée d'un tel morphisme est équivalente à celle de $d$ sections globales de $p^{-1}(\eta)$. Soit $(T_i)_{i\in I}$ un recouvrement ouvert de $T$ par des ouverts affines $T_i=\Spec B_i$, où $B_i$ est une $k$-algèbre de type fini. Comme la donnée d'une section globale de $p^{-1}(\eta)$ est équivalente à celle de $|I|$ sections de $p^{-1}(\eta)$ au-dessus des $T_i\times_{k}\Spec \kappa(\eta)$ avec conditions de recollement, il nous suffit de prouver l'existence d'un tel morphisme pour $T$ affine. Soient $B$ la $k$-algèbre des sections globales de $T$ et $K$ le corps des fractions de $A$. Soit $(s_i)_{1\leq i\leq d}$ $d$ éléments de $B_K=B\otimes_A K$. Il existe alors $g\in A$ tel que $(s_i)_{1\leq i\leq d}\in B_g=B\otimes_{A}A_g$.

\textit{Prouvons que, quitte à restreindre $U$, le $U$-morphisme $p^{-1}(U)\rightarrow \mathbf{A}_U$ est encore une immersion fermée}. On peut supposer que $T=\Spec B$. Si le $K$-morphisme d'algèbres $K[X_1,\ldots,X_d]\rightarrow B_K$ est surjectif, il existe un élément $g'\in A$ tel que la restriction du $U$-morphisme $p^{-1}(U)\rightarrow \mathbf{A}_U$ à $T\times_k D(gg')$ soit déjà une immersion fermée, en appliquant le raisonnement précédent aux générateurs de $B$ comme $A$-algèbre de type fini. En particulier, ceci entraine que le $k$-schéma $T$ est affine. On note $T=\Spec B$, où $B$ est une $k$-algèbre de type fini.

\textit{Montrons que, quitte à restreindre $U$, le noyau du morphisme d'algèbres $A_{gg'}\rightarrow B\otimes_A A_{gg'}$ est nul}. On applique le raisonnement précédent aux générateurs du noyau, pour trouver un élément $g''\in A$ tel que la restriction de notre morphisme de schémas $T\times_k D(gg'g'')$ soit un isomorphisme.

\end{proof}

\begin{proof}[Démonstration de la proposition \ref{proprietes de la mesure}]

Soit $A$ un cylindre de $\G(\!X)$ de rang $m$. Le lemme clé \ref{2} assure que pour tout $e\in\mathbf{N}$ l'ensemble $A\cap\G^{(e)}(\!X)$ est un cylindre stable de rang $\max(m, c_\!Xe)$. Il nous faut donc montrer que la suite $(\mu_{0,\!X}(A\cap\G^{(e)}(\!X)))_{e\in\mathbf{N}}$ converge dans $\widehat{\!M}$, \textit{i.e.} est de Cauchy dans $\widehat{\!M}$. Posons $A^{(e)}:=A\cap\G^{(e)}(\!X)$ et $A_{(e)}:=\G(\!X)\backslash (A\cap\G^{(e)}(\!X))=A\cap \pi_{e,\!X}^{-1}(\G_e(\!X_{\mathrm{sing}}))$. Soient $\varepsilon>0$ et $q=E(-log(\varepsilon)/log 2)+1$. On veut montrer qu'il existe $e_0\in\mathbf{N}$ tel que pour tous $e'\geq e\geq e_0$
$$
[\pi_{n,\!X}(A^{(e)}]\mathbf{L}^{-(n+1)d}-[\pi_{n,\!X}(A^{(e')}]\mathbf{L}^{-(n+1)d}\in F^{q+1}\widehat{\!M}.
$$
Par additivité de $[\hskip 1mm]$ sur les constructibles, on obtient, pour $n\geq\mathrm{max}(m,c_\!Xe,c_\!Xe')$:
$$
[\pi_{n,\!X}(A^{(e)})]\mathbf{L}^{-(n+1)d}=[\pi_{n,\!X}(A)]\mathbf{L}^{-(n+1)d} -[\pi_{n,\!X}(A_{(e)})]\mathbf{L}^{-(n+1)d}.
$$
En particulier, ceci entra\^\i ne que
$$
[\pi_{n,\!X}(A^{(e)}]\mathbf{L}^{-(n+1)d}-[\pi_{n,\!X}(A^{(e')}]\mathbf{L}^{-(n+1)d}
$$
$$
=([\pi_{n,\!X}(A^{(e)})\cap \pi_{n,\!X}(A_{(e)})]-[\pi_{n,\!X}(A^{(e')})\cap \pi_{n,\!X}(A_{(e')})])\mathbf{L}^{-(n+1)d}
$$
$$
+ ([\pi_{n,\!X}(A_{(e')})]-[\pi_{n,\!X}(A_{(e)})])\mathbf{L}^{-(n+1)d}.
$$
Si $e'\geq e\geq \gamma(q)$, où $\gamma$ désigne la fonction de Greenberg pour $\!X_{\mathrm{sing}}$, le lemme \ref{3} assure que 
$$
[\pi_{n,\!X}(A_{(e')})]\mathbf{L}^{-(n+1)d}\in F^{q+1}\widehat{\!M}
$$
et
$$
[\pi_{n,\!X}(A_{(e)})]\mathbf{L}^{-(n+1)d}\in F^{q+1}\widehat{\!M}
$$
Ceci implique que la suite $([\pi_{n,\!X}(A^{(e)})]\mathbf{L}^{-(n+1)d})_{e\in\mathbf{N}}$ est de Cauchy dans $\widehat{\!M}$.

Les autres assertions sont vérifiées par $\mu_{0,\!X}$ et sont stables par passage à la limite.
\end{proof}

\begin{lem}
Soit $A$ un cylindre de $\G(\!X)$ de rang $m$. La suite 
$$
([\pi_{n,\!X}(A)]\mathbf{L}^{-(n+1)d})_{n\in\mathbf{N}}
$$
converge dans $\widehat{\!M}$ vers $\tilde{\mu}_\!X(A)$.
\end{lem}

\begin{proof}
Soit $e\geq 0$, posons $A^{(e)}:=A\backslash \pi_{e,\!X}^{-1}(\G_e(\!X_{\mathrm{sing}}))$. Soit $\gamma$ la fonction de Greenberg pour $\!X_{\mathrm{sing}}$. Le lemme \ref{3} assure que, pour tout $n\geq\mathrm{max}(\gamma(e),m,c_\!Xe)$,
$$
[\pi_{n,\!X}(A)]\mathbf{L}^{-(n+1)d}-[\pi_{n,\!X}(A^{(\gamma(e))})]\mathbf{L}^{-(n+1)d}
$$
appartient à $F^{e+1}\!M_{\mathrm{loc}}$. Par définition, $A^{(e)}\subset \G^{(e)}(\!X)$. Il découle du lemme \ref{2} que $A^{(\gamma(e))}$ est stable de rang $\mathrm{max}(m,c_\!X\gamma(e))$. Par conséquent, pour tous $n$ et $n'\geq \mathrm{max}(m,c_\!X\gamma(e))$,
$$
[\pi_{n,\!X}(A^{(\gamma(e))})]\mathbf{L}^{-(n+1)d}=[\pi_{n',\!X}(A^{(\gamma(e))})]\mathbf{L}^{-(n'+1)d}.
$$
On en déduit que pour tous $n$ et $n'\geq \mathrm{max}(\gamma(e),m,c_\!Xe)$, $[\pi_{n,\!X}(A)]\mathbf{L}^{-(n+1)d}-[\pi_{n',\!X}(A)]\mathbf{L}^{-(n'+1)d}
$ appartient à $F^{e+1}\!M_{\mathrm{loc}}$. Ceci signifie que la suite 
$$
([\pi_{n,\!X}(A)]\mathbf{L}^{-(n+1)d})_{n\in\mathbf{N}}
$$
est une suite de Cauchy et donc converge vers $\!L$ dans $\widehat{\!M}$. Comme, par définition, $\tilde{\mu}_\!X(A^{(\gamma(e))})=[\pi_{n,\!X}(A^{(\gamma(e))})]\mathbf{L}^{-(n+1)d}$ pour $n\geq \mathrm{max}(m,c_\!X\gamma(e))$, on en conclut que la suite de terme général $\tilde{\mu}_\!X(A^{(\gamma(e))})$ converge vers $\!L$ dans $\widehat{\!M}$. Comme, par définition, $(\tilde{\mu}_\!X(A^{(\gamma(e))}))$ converge vers $\tilde{\mu}_\!X(A)$, on en déduit que $\mu_\!X(A)=\!L$.
\end{proof}


\section{Ensembles mesurables}


Supposons désormais que $\!X$ est un $R$-schéma formel admissible et réduit.


\subsection{La définition.}


\begin{defi}
Soit $A$ un sous-ensemble de $\G(\!X)$. On dit que $A$ est mesurable dans $\G(\!X)$ si, pour tout $\varepsilon>0$, il existe un ensemble $I(\varepsilon)$ au plus dénombrable et une famille de cylindres $(A_i(\varepsilon))_{i\in I\cup\{0\}}$ tels que
\begin{enumerate}

\item  $A\triangle A_0(\varepsilon)\subset\displaystyle\bigcup_{i\in I}A_i(\varepsilon)$, où $\triangle$ désigne la différence symétrique. 

\item  $\lVert\mu(A_i(\varepsilon))\rVert< \varepsilon$, pour tout $i\in I$.
\end{enumerate}
Dans ce cas, on dit qu'une telle famille de cylindres $(A_i(\varepsilon))_{i\in I(\varepsilon)\cup\{0\}}$ est une $\varepsilon$-approximation cylindrique de $A$. Le cylindre $A_0(\varepsilon)$ est appelé la partie principale de l'approximation cylindrique $A_{I(\varepsilon)}$. On dit que $A$ est fortement mesurable si de plus on peut choisir $A_0(\varepsilon)\subset A$, pour tout $\varepsilon>0$. 
\end{defi}

\begin{exes}
\label{pro-cylindre=mesurable}

\begin{enumerate}
\item  Tout cylindre $C$ de $\G(\!X)$ est mesurable (même fortement mesurable). En effet, si $\varepsilon>0$, la famille de cylindres $(C,\emptyset,\emptyset,\emptyset,\ldots)$ est une $\varepsilon$-approximation cylindrique de $C$.

\item Soit $\!Z\hookrightarrow\!X$ un sous-$R$-schéma formel fermé d'un $R$-schéma sttf lisse sur $R$ de pure dimension $d$. Alors $\G(\!Z)$ est un ensemble mesurable de $\G(\!X)$. En effet, si $\!Z$ ne contient aucune composante de $\!X$, le résultat est une conséquence du lemme clé \ref{3}. Sinon $\!Z$ s'écrit comme
$$
\!Z=(\sqcup_{i\in I} \!X_i)\sqcup(\tilde{\!Z})
$$
où les $\!X_i$ sont des composantes connexes de $\!X$ et $\tilde{\!Z}$ un sous-$R$-schéma formel fermé de $\!X$ de dimension au plus $d-1$. On a alors
$$
\G(\!Z)\triangle\sqcup(\sqcup_{i\in I}\G(\!X_i))\subset \pi_{n,\!X}^{-1}(\G_n(\tilde{Z}))
$$
pour tout $n\in\mathbf{N}$. En effet, soit $z\in\G(\!Z)$ tel que $z\not\in\sqcup_{i\in I} \!X_i$ et $\pi_{n,\!X}(z)\in\G_n(\!X_{i_0})$. Alors $\pi_{0,\!X}(z)\in\G_n(\!X_{i_0})$ et, comme $\!X_{i_0}$ est ouvert, le point $z$ appartient à $G_n(\!X_{i_0})$.Le résultat découle alors du lemme \ref{3}. Soient $\varepsilon>0$ et $q=E(-log(\varepsilon)/log(2))+1$. Pour tout $n\geq \gamma_\!Z(q)$, on a
$$
[\pi_{n,\!X}(\pi_{n,\!X}^{-1}(\G_n(\tilde{\!Z})))]\mathbf{L}^{-(n+1)d}\in F^{q+1}\widehat{\!M}.
$$
\end{enumerate}
\end{exes}

On d\'esignera par $\mathbf{D}_\!X$ (\textit{resp.} $\mathbf{D}_\!X^{f}$) l'ensemble des ensembles mesurables (\textit{resp.} fortement mesurables) de $\G(\!X)$.

\begin{prop}
\label{mesurables sont un anneau booleen}
Les ensembles de parties $\mathbf{D}_\!X$ et $\mathbf{D}_\!X^{f}$  v\'erifient les propriétés suivantes:

\begin{enumerate}

\item $\mathbf{D}_\!X^{f}\subset\mathbf{D}_\!X$.

\item Les ensembles $\G(\!X)\in\mathbf{D}_\!X^f$ et $\emptyset\in\mathbf{D}_\!X^f$.

\item L'ensemble $\mathbf{D}_\!X$ (resp. $\mathbf{D}_\!X^f$) est stable par intersection finie.

\item L'ensemble $\mathbf{D}_\!X$ (resp. $\mathbf{D}_\!X^f$) est stable par r\'eunion finie.

\item L'ensemble $\mathbf{D}_\!X$ est stable par passage au compl\'ementaire.

\item Si $A$ et $B$ sont deux parties mesurables de $G(\!X)$, alors $A\backslash B$ est encore une partie mesurable de $\G(\!X)$.

\end{enumerate}
En particulier, $\mathbf{D}_\!X$ est un anneau bool\'een.
\end{prop}

\begin{proof}

Les ensembles $\G(\!X)$ et $\emptyset$ sont des cylindres, donc des ensembles fortement mesurables. Soit $\varepsilon>0$. Soit $A^{(i)}$ pour $1\leq i\leq n$ une famille finie d'ensembles mesurables de $\G(\!X)$. Pour tout $1\leq i\leq n$, soit $(A_j^{(i)}(\varepsilon))_{j\in I_i\cup \{0\}}$ une $\varepsilon$-approximation cylindrique de $A^{(i)}$. La famille de cylindres définie par $A_0(\varepsilon):=\cup_{i=1}^{n} A^{(i)}_0(\varepsilon)$, et les $A^{(i)}_j(\varepsilon)$ pour tout $1\leq i\leq n$ et tout $j\in I_i$, est une $\varepsilon$-approximation cylindrique de $A:=\cup_{i=1}^{n} A^{(i)}$. De même, la famille de cylindres définie par $B_0(\varepsilon):=\cap_{i=1}^{n} A^{(i)}_0(\varepsilon)$, et les $A^{(i)}_j(\varepsilon)$ pour tout $1\leq i\leq n$ et tout $j\in I_i$, est une $\varepsilon$-approximation cylindrique de $B:=\cap_{i=1}^{n} A^{(i)}$. Le point 5 découle de l'égalité suivante
$$
A\triangle A_0(\varepsilon)=(\G(\!X)\backslash A)\triangle (\G(\!X)\backslash A_0(\varepsilon)).
$$ 
Le point 6 est une conséquence formelle des points 3 et 5 et la suite de cylindres définie par $C_0(\varepsilon):=A^{(1)}_0(\varepsilon)\backslash A^{(2)}_0(\varepsilon)$, les $A^{(1)}_j(\varepsilon)$ pour tout $j\in I_1$ et les $A^{(2)}_i(\varepsilon)$ pour tout $i\in I_2$ est une $\varepsilon$-approximation cylindrique de $C:=A^{(1)}\backslash A^{(2)}$.
\end{proof}


\subsection{La mesure motivique pour les ensembles mesurables.}


\begin{thm}
\label{theoreme de Batyrev}
Soit $\!X$ un $R$-schéma formel admissible et réduit.
\begin{enumerate}

\item Pour tout ensemble mesurable $A\in\mathbf{D}_\!X$, la limite
$$
\mu_\!X(A):={\displaystyle{\lim_{\varepsilon\rightarrow 0}}}\hskip0.5mm\tilde{\mu}_\!X(A_0(\varepsilon))
$$
où $A_0(\varepsilon)$ désigne la partie principale d'une $\varepsilon$-approximation cylindrique, existe dans $\widehat{\!M}$. En outre, cette limite est ind\'ependante du choix des approximations cylindriques de $A$.

\item Si $A\in\mathbf{C}_{\!X}$, alors ${\mu}_\!X(A)$ et $\tilde{\mu}_\!X(A)$ co\"{i}ncident dans $\widehat{\!M}$.

\item L'application $A\mapsto{\mu}_\!X(A)$ est $\sigma$-additive sur $A\in\mathbf{D}_\!X$.

\item Pour $A$ et $B$ dans $A\in\mathbf{D}_\!X$, $\lVert{\mu}_\!X(A\cup B)\rVert\leq \max(\lVert{\mu}_\!X(A)\rVert,\lVert{\mu}_\!X(B)\rVert)$. Si $A\subset B$, $\lVert{\mu}_\!X(A)\rVert\leq\lVert{\mu}_\!X(B)\rVert$.

\end{enumerate}

\end{thm}

\begin{proof}
On omettra l'indice $\!X$ dans la preuve. La preuve du point 1 est l'adaptation d'une preuve de Batyrev (\textit{cf} \cite{ba}), établie pour démontrer un résultat analogue. Soit $A$ un ensemble mesurable. Soient  $\varepsilon>0$ et $\varepsilon'>0$ deux nombres r\'eels. Soient $A_{I(\varepsilon)}$ et $A'_{I'(\varepsilon')}$ deux approximations cylindriques de $A$.

Par d\'efinition, on a la relation 
$$
A_0(\varepsilon)\triangle A'_0(\varepsilon')\subset\Bigl(\displaystyle\bigcup_{i\in I(\varepsilon)}A_i(\varepsilon)\Bigr)\cup\Bigl(\displaystyle\bigcup_{i\in I'(\varepsilon')}A'_i(\varepsilon')\Bigr).
$$ 
Comme $A_0(\varepsilon)\triangle A'_0(\varepsilon')$ est un cylindre, le lemme \ref{quasi-compacite de la topologie constructible} assure qu'il existe deux entiers $L(\varepsilon)$ et $L(\varepsilon')$ tels que 
$$
A_0(\varepsilon)\triangle A'_0(\varepsilon')\subset\Bigl(\displaystyle\bigcup_{i= 1}^{L(\varepsilon)}A_i(\varepsilon)\Bigr)\cup\Bigl(\displaystyle\bigcup_{i= 1}^{L(\varepsilon')}A'_i(\varepsilon')\Bigr).
$$ 
Les propri\'et\'es de la norme sur $\widehat{\!M}$ donnent la relation 
$$
\lVert\mu(A_0(\varepsilon)\triangle A'_0(\varepsilon'))\rVert\leq\mathrm{max}(\varepsilon,\varepsilon').
$$ 
En utilisant les inclusions 
$$
A_0(\varepsilon)\backslash(A_0(\varepsilon)\cap A'_0(\varepsilon'))\subset A_0(\varepsilon')\triangle A'_0(\varepsilon')
$$ 
$$
A'_0(\varepsilon')\backslash(A_0(\varepsilon)\cap A'_0(\varepsilon'))\subset A_0(\varepsilon)\triangle A'_0(\varepsilon')
$$ 
on obtient les in\'egalit\'es 
$$
\lVert\mu(A_0(\varepsilon)\backslash(A_0(\varepsilon)\cap A'_0(\varepsilon')))\rVert\leq\mathrm{max}(\varepsilon,\varepsilon')
$$ 
$$
\lVert\mu(A'_0(\varepsilon)\backslash(A_0(\varepsilon)\cap A'_0(\varepsilon')))\rvert\leq\mathrm{max}(\varepsilon,\varepsilon').
$$ 
En utilisant la proposition \ref{proprietes de la mesure} et les d\'ecompositions $A_0(\varepsilon)=(A_0(\varepsilon)\backslash(A_0(\varepsilon)\cap A'_0(\varepsilon')))\sqcup (A_0(\varepsilon)\cap A'_0(\varepsilon'))$, $A'_0(\varepsilon)=(A'_0(\varepsilon)\backslash(A_0(\varepsilon)\cap A'_0(\varepsilon')))\sqcup (A_0(\varepsilon)\cap A'_0(\varepsilon'))$, on obtient 
$$
\lVert\mu(A_0(\varepsilon))-\mu(A'_0(\varepsilon'))\rVert\leq\mathrm{max}(\varepsilon,\varepsilon').
$$ 
Ce r\'esultat suffit pour montrer le th\'eor\`eme. En effet, si $A'_{I'(\varepsilon')}=A_{I(\varepsilon')}$, on vient de montrer que la famille $(A_{I(\varepsilon)})_{\varepsilon>0}$ est un filtre de Cauchy, donc converge dans $\widehat{\!M}$. En outre, soient $a$ et $a'$ les limites respectives des suites $(\mu(A_0(1/n))_{n\in\mathbf{N}}$ et $(\mu(A'_0(1/n))_{n\in\mathbf{N}}$. On a l'inégalité 
$$
\lVert a-a'\rVert\leq \lVert a-\mu(A_0(1/n))\rVert+\lVert \mu(A_0(1/n))-\mu(A'_0(1/m))\lVert+\lVert a'-\mu(A_0(1/m))\lVert.
$$
Il est alors clair que $a=a'$. Par cons\'equent, les deux limites co\"\i ncident. 

Soit $C$ un cylindre de $\G(\!X)$. Le point 2 découle du fait que, pour tout $\varepsilon>0$, on peut choisir une $\varepsilon$-approximation de $C$ de partie principale égale à $C$.

Les points 3 et 4 sont démontrés ci-dessous. Pour le point 3, \textit{cf} corollaire \ref{description des reunions de mesurables}/1. Soient $A$ et $B$ deux parties mesurables de $\G(\!X)$. Soit $\varepsilon>0$. Il existe une $\varepsilon$-approximation cylindrique de $A\cup B$ de partie principale $A_0(\varepsilon)\cup B_0(\varepsilon)$, où $A_0(\varepsilon)$ (\textit{resp.} $B_0(\varepsilon)$) est la partie principale d'une $\varepsilon$-approximation de $A$ (\textit{resp.} $B$). La première partie de 4 en découle. La dernière assertion de 4 est démontrée dans la preuve du lemme \ref{pt 4}.

\end{proof}


\subsection{Les propri\'et\'es g\'en\'erales}
\label{proprietes de la mesure sur les mesurables}


\subsubsection{Les propriétés des approximations cylindriques}

\begin{lem}
\label{approximation disjointe}
Soient $A$ et $B$ deux ensembles mesurables disjoints, alors, pour tout $\varepsilon>0$, il existe une $\varepsilon$-approximation cylindrique $A_\ast(\varepsilon)$ (resp. $B_\ast(\varepsilon)$) de $A$ (resp. $B$) telle que $A_0(\varepsilon)\cap B_0(\varepsilon)=\emptyset$.
\end{lem}

\begin{proof}
Soit $\varepsilon>0$. Par définition, il existe une $\varepsilon$-approximation cylindrique $(A'_i(\varepsilon))_{i\in I\cup\{0\}}$ (resp. $(B'_i(\varepsilon))_{i\in J\cup\{0\}}$) de $A$ (resp. $B$). La suite de cylindres de $\G(\!X)$ définie par $A_0(\varepsilon):=A'_0(\varepsilon)\backslash B'_0(\varepsilon)$, par les $A_i(\varepsilon)$ pour tout $i\in I$ et par les $B'_j(\varepsilon)$ pour tout $j\in J$ (\textit{resp.} par $B'_j(\varepsilon)$ pour tout $j\in J\cup\{0\}$) est une $\varepsilon$-approximation cylindrique de $A$ (\textit{resp.} de $B$). Ces deux suites répondent à la question.
\end{proof}

\begin{cor}
Soit $(A^{(i)})_{1\leq i\leq n}$ une famille finie de parties mesurables deux \`a deux disjointes, alors, pour tout $1\leq i\leq n$, pour tout $\varepsilon>0$, il existe une $\varepsilon$-approximation cylindrique $A_{I_i}^{(i)}$ de $A^{(i)}$, telle que si $i\not= j$, $A_0^{(i)}(\varepsilon)\cap A_0^{(j)}(\varepsilon)=\emptyset$.
\end{cor}

\begin{proof}
Soit $\varepsilon>0$. L'assertion se prouve par récurrence sur $n$. Supposons que ce cardinal soit égal à $n+1$. Considérons $n+1$ parties mesurables de $\G(\!X)$ deux à deux disjointes. Par hypothèse de récurrence, on peut trouver, pour tout $1\leq i\leq n$, des $\varepsilon$-approximations cylindriques $A^{(i)}_{I_i}$ de $A^{(i)}$ de parties principales deux à deux disjointes. La proposition \ref{mesurables sont un anneau booleen}/4 assure qu'il existe une $\varepsilon$-approximation cylindrique de $\cup_{1\leq i\leq n} A^{(i)}$ de partie principale $\cup_{1\leq i\leq} A^{(i)}_0(\varepsilon)$. La preuve du lemme \ref{approximation disjointe} construit alors une $\varepsilon$-approximation cylindrique de $A^{(n+1)}$ de partie principale d'intersection vide avec $\cup_{1\leq i\leq} A^{(i)}_0(\varepsilon)$.
\end{proof}

\begin{lem}
\label{approximation contenue}
Si $A$ et $B$ sont deux ensembles mesurables tels que $A\subset B$, alors, pour tout $\varepsilon>0$, il existe une $\varepsilon$-approximation cylindrique $A_\ast$ (\textit{resp.} $B_\ast$) de $A$ (\textit{resp.} de $B$) telle que $A_0(\varepsilon)\subset B_0(\varepsilon)$. En particulier, si $B$ est fortement mesurable, on peut supposer que $A_0(\varepsilon)\subset B_0(\varepsilon)\subset B$.
\end{lem}

\begin{proof}
Soit $\varepsilon>0$. Soit $A'_\ast(\varepsilon)$ (\textit{resp.} $B'_\ast(\varepsilon)$) une $\varepsilon$-approximation de $A$ (resp. $B$). La suite de cylindres $A_\ast(\varepsilon)$ (\textit{resp.} $B_\ast(\varepsilon)$) définie par $A_0(\varepsilon):=A'_0(\varepsilon)\cap B'_0(\varepsilon)$, par les $A'_j(\varepsilon)$ pour tout $j\not=0$ et les $B'_j(\varepsilon)$ pour tout $j\not=0$ (\textit{resp.} par les $B'_j(\varepsilon)$) est une $\varepsilon$-approximation cylindrique de $A$ (resp. $B$) vérifiant la propriété de l'énoncé.
\end{proof}

\subsubsection{Les propriétés de la mesure sur les ensembles mesurables.}

\begin{lem}
\label{pt 4}
Si $A$ et $B$ sont deux ensembles mesurables de $\G(\!X)$ et si $A\subset B$, alors on a $\lVert\mu_\!X(A)\rVert\leq\lVert\mu_\!X(B)\rVert$. 
\end{lem}

\begin{proof}
D'apr\`es le lemme \ref{approximation contenue}, on peut trouver une $\varepsilon$-approximation cylindrique de $A$ et une $\varepsilon$-approximation de $B$ telles que $A_0(\varepsilon)\subset B_0(\varepsilon)$. Le r\'esultat d\'ecoule alors du fait que $\lVert\mu_\!X(A_0(\varepsilon))\rVert\leq\lVert\mu_\!X(B_0(\varepsilon))\rVert$ et de la continuit\'e de la norme.
\end{proof}

\begin{prop}
\label{mesurable comme limite}
Soit $(A^{(i)})_{i\in\mathbf{N}}$ une suite d'ensembles mesurables (resp. fortement mesurables) tels que la suite $(\mu_\!X(A^{(i)}))_{i\in\mathbf{N}}$ tende vers 0 dans $\widehat{\!M}$. Alors ${\cup_{i\in\mathbf{N}}}A^{(i)}$ est mesurable (resp. fortement mesurable).
\end{prop}

\begin{proof}
On posera $\mu$ pour $\mu_\!X$ dans cette preuve. Soit $\varepsilon>0$. Il s'agit de construire une $\varepsilon$-approximation cylindrique de $\cup_{i\in\mathbf{N}} A^{(i)}$. Pour tout $\varepsilon'>0$, pour tout $i\in\mathbf{N}$, soit $A^{(i)}_{I_i(\varepsilon')}$ une $\varepsilon'$-approximation cylindrique de $A^{(i)}$. On a montr\'e dans la preuve du th\'eor\`eme \ref{theoreme de Batyrev} que: 
$$
\lVert\mu(A_0^{(i)}(\varepsilon/2))-\mu(A_0^{(i)}(1/n))\rVert\leq \max(\varepsilon/2, 1/n).
$$ 
En particulier, ceci entraine que, pour tout $i\in\mathbf{N}$,
$$
\lVert\mu(A_0^{(i)}(\varepsilon/2))\rVert\leq \varepsilon/2 +\lVert\mu(A^{(i)})\rVert.
$$
Enfin, par hypoth\`ese, il existe $i_0\in\mathbf{N}$ tel que, pour tout $i> i_0$, $\lVert\mu(A^{(i)})\rVert<\varepsilon/2$. Donc 
$$
\lVert\mu(A_0^{(i)}(\varepsilon/2))\rVert<\varepsilon.
$$ 
La suite de cylindres définie par $B_0(\varepsilon):=\cup_{i=0}^{i_0}A^{(i)}_0(\varepsilon/2)$, par les $A_{j}^{(i)}(\varepsilon)$ pour $i\in\mathbf{N}$ et $j\in I_i$, et les $A^{(i+i_0+1)}_0(\varepsilon/2)$ pour tout $i\in\mathbf{N}$, est une $\varepsilon$-approximation cylindrique de $B:=\cup A^{(i)}$, grâce au choix de $i_0$.

\end{proof}

\begin{prop}
\label{mesurable comme limite bis}
Soit $A:=\cup_{i\in\mathbf{N}}A^{(i)}$ une r\'eunion d\'enombrable d'ensembles mesurables de $\G(\!X)$ deux \`a deux disjoints . Si on suppose que $A$ est un ensemble mesurable de $\G(\!X)$, alors $\mu_\!X(A)=\Sigma_{i\in \mathbf{N}} \mu_\!X(A^{(i)})$.
\end{prop}

\begin{proof}
Pour tout $\delta>0$, $A_{I(\delta)}(\delta)$ désigne une $\delta$-approximation cylindrique de $A$ et $A^{(j)}_{I_j(\delta)}$ une $\delta$-approximation cylindrique de $A^{(j)}$ pour tout $j\in\mathbf{N}$.

Soit $1>\varepsilon>0$. Par définition de la mesure $\mu_\!X$, il existe $\varepsilon>\varepsilon'>0$ tel que $\lVert\mu_\!X(A)-\mu_\!X(A_0(\varepsilon'))\rVert<\varepsilon$. Comme $A\triangle A_0(\varepsilon')\subset\cup_{i\in I(\varepsilon')}A_i(\varepsilon')$, 
$$
A_0(\varepsilon')\subset (\cup_{j\in\mathbf{N}}A^{(j)})\cup(\cup_{i\in I(\varepsilon')} A_i(\varepsilon')).
$$
De même, comme $A^{(j)}\triangle A_0^{(j)}(\varepsilon')\subset\cup_{i\in I_j(\varepsilon')}A_i^{(j)}(\varepsilon')$,
$$
A_0(\varepsilon')\subset (\cup_{j\in\mathbf{N},i\in I_j(\varepsilon')\cup\{0\}}A_i^{(j)})\cup(\cup_{i\geq 1} A_i(\varepsilon')).
$$
Le lemme \ref{quasi-compacite de la topologie constructible} assure l'existence d'un entier $m(\varepsilon')\in\mathbf{N}$ tel que pour tout $n\geq m(\varepsilon')$
$$
A_0(\varepsilon')\subset (\cup_{j\leq n}A_0^{(j)})\cup(\cup_{i\in M} C_i)
$$
où $M$ est un ensemble fini et les $C_i$ sont des cylindres de $\G(\!X)$ tels que $\lVert\mu_\!X(C_i))\rVert<\varepsilon'$ pour tout $i\in M$. En remarquant que la suite de cylindres de $\G(\!X)$ définie par $B_0^{(j)}(\varepsilon'):=A_0^{(j)}(\varepsilon')\cap A_0(\varepsilon')$, par les $A_i(\varepsilon')$ pour tout $i\in I$ et les $A_i^{(j)}(\varepsilon')$ pour tout $i\in I_j$, est une $\varepsilon'$-approximation cylindrique de $A^{(j)}$ pour tout $j\in\mathbf{N}$, et que, pour tout $0\leq j\leq n$, l'on peut choisir les $A_0^{(j)}(\varepsilon')$ disjoints deux à deux, on a dans $\widehat{\!M}$ l'égalité:
$$
\mu_\!X(A_0(\varepsilon'))=\sum_{j=0}^n\mu_\!X(B_0^{(j)}(\varepsilon'))+\sum_{i\in I}\mu_\!X(C'_i)
$$
où les cylindres $C'_i:=C_i\cap A_0(\varepsilon')$ de $\G(\!X)$ sont tels que  $\lVert\mu_\!X(C'_i))\rVert<\varepsilon'$. En particulier, ceci entra\^\i ne que pour tout $n\geq m(\varepsilon')$:
$$
\mu_\!X(A_0(\varepsilon'))-\sum_{j=0}^n\mu_\!X(B_0^{(j)}(\varepsilon'))\in F^q\widehat{\!M}
$$
avec $q=E(-log(\varepsilon)/log 2)$ et donc, grâce au choix de $\varepsilon'$, que
$$
\mu_\!X(A)-\sum_{j=0}^n\mu_\!X(B_0^{(j)}(\varepsilon'))\in F^q\widehat{\!M}.
$$
Par ailleurs, par définition de la mesure, il existe $\varepsilon'>\varepsilon''>0$ tel que
$$
\sum_{j=0}^n\mu_\!X(B_0^{(j)}(\varepsilon''))-\sum_{j=0}^n\mu_\!X(A^{(j)})\in F^q\widehat{\!M}.
$$

En outre, par un argument utilisé dans la preuve du théorème \ref{theoreme de Batyrev},
$$
\sum_{j=0}^n\mu_\!X(B_0^{(j)}(\varepsilon'))-\sum_{j=0}^n\mu_\!X(B_0^{(j)}(\varepsilon''))\in F^q\widehat{\!M}.
$$
En additionnant ces trois dernières égalités, on obtient que pour tout $n\geq m(\varepsilon')$:
$$
\lVert\mu_\!X(A)-\sum_{j=0}^n\mu_\!X(A^{(j)})\lVert<\varepsilon.
$$
\end{proof}

\begin{cor}
\label{reunion et intersection de mesurables}
\label{description des reunions de mesurables}

\begin{enumerate}

\item  Soit $(A_i)_{i\in \mathbf{N}}$ une famille de mesurables deux \`a deux disjoints, alors leur r\'eunion $\cup_{i\in\mathbf{N}} A_i$ est mesurable si et seulement si la suite des volumes $(\mu_\!X(A_i))_{i\in\mathbf{N}}$ tend vers 0 quand $i$ tend vers l'infini. Dans ce cas, on a la relation:$$\mu_{\!X}(\cup_{n\in\mathbf{N}}A_n)=\sum_{n\in\mathbf{N}}\mu_\!X(A_n).$$

\item Soit $A$ un sous-ensemble de $\G(\!X)$ tel qu'il existe une suite d'ensembles mesurables $(A_n)_{n\in\mathbf{N}}$ de $\G(\!X)$ vérifiant les conditions suivantes:
\begin{enumerate}

\item pour tout $n\geq 0$, $A_{n+1}\subset A_n$.

\item $A=\cap_{n\in\mathbf{N}} A_n$.

\end{enumerate}

Alors $A$ est mesurable si et seulement si la suite $(\mu_\!X(A_n))_{n\in\mathbf{N}}$ converge dans $\widehat{\!M}$. Dans ce cas, on a 
$$
\mu_\!X(A)=\lim_{n\rightarrow+\infty}\mu_\!X(A_n).
$$

\item Soit $B$ un ensemble mesurable de $\G(\!X)$ tel qu'il existe une suite d'ensembles mesurables $(B_n)_{n\in\mathbf{N}}$ de $\G(\!X)$ vérifiant les conditions suivantes:
\begin{enumerate}

\item pour tout $n\geq 0$, $B_{n}\subset B_{n+1}$.

\item $B=\cup_{n\in\mathbf{N}} B_n$.

\end{enumerate}

Alors $B$ est mesurable si et seulement si la suite $(\mu_\!X(B_n))_{n\in\mathbf{N}}$ converge dans $\widehat{\!M}$. Dans ce cas, on a 
$$
\mu_\!X(B)=\lim_{n\rightarrow+\infty}\mu_\!X(B_n).
$$ 

\end{enumerate}
\end{cor}

\begin{proof}
Dans cette preuve, on omettra les indices $\!X$. Le point \textbf{1}. est une simple retraduction des propositions \ref{mesurable comme limite} et \ref{mesurable comme limite bis} et du fait que la topologie de $\widehat{\!M}$ est ultramétrique. \textbf{2}. Pour tout $n\in\mathbf{N}$, posons $C_n=A_n\backslash A_{n+1}$. Les $C_n$ sont alors des ensembles mesurables disjoints. Par ailleurs,
$$
\cup_{n\in\mathbf{N}} C_n=A_0\backslash(\cap_{m\geq 0} A_m).
$$
Cette égalité entraine en particulier que $A$ est mesurable si et seulement si la réunion $\cup_{n\in\mathbf{N}} C_n$ l'est. Or le point 1. implique que $
\cup_{n\in\mathbf{N}} C_n$ est mesurable si et seulement si la suite $(\mu(C_n))$ tend vers 0 dans $\widehat{\!M}$. Par additivité de la mesure, on a l'égalité $\mu(C_n)=\mu(A_n)-\mu(A_{n+1})$. Comme la topologie de $\widehat{\!M}$ est ultramétrique, la suite $(\mu(C_n))$ tend vers 0 dans $\widehat{\!M}$ si et seulement si la suite $(\mu(A_n))$ est de Cauchy dans $\widehat{\!M}$. Enfin, comme $\widehat{\!M}$ est complet, $\cup_{n\in\mathbf{N}} C_n$ est mesurable si et seulement si la suite $(\mu(A_n))$ converge dans $\widehat{\!M}$. Dans ce cas, le lemme \ref{reunion et intersection de mesurables} assure que $\mu(A)=\lim \mu(A_n)$. \textbf{3}. Posons $C_0=B_0$, $C_n=B_n-B_{n-1}$, pour $n\geq 1$. Les $C_n$ appartiennent \`a $\mathbf{D}_\!X$ et sont disjoints deux \`a deux. En outre, comme $B_n=\sqcup_{i=1}^n B_i$ on a $B=\sqcup_{n\in\mathbf{N}} C_n$. L'assertion 1. assure que $B$ est mesurable si et seulement si la suite $(\mu(C_n))$ tend vers 0 dans $\widehat{\!M}$. Par additivité de la mesure, on a l'égalité $\mu(C_n)=\mu(B_n)-\mu(B_{n-1})$ pour $n\geq 1$. Comme la topologie de $\widehat{\!M}$ est ultramétrique, la suite $(\mu(C_n))$ tend vers 0 dans $\widehat{\!M}$ si et seulement si la suite $(\mu(A_n))$ est de Cauchy dans $\widehat{\!M}$. Enfin, comme $\widehat{\!M}$ est complet, $\cup_{n\in\mathbf{N}} C_n$ est mesurable si et seulement si la suite $(\mu(B_n))$ converge dans $\widehat{\!M}$. Dans ce cas, le lemme \ref{reunion et intersection de mesurables} assure que $\mu(B)=\lim \mu(B_n)$.
\end{proof}


\subsection{Les ensembles de mesure nulle}
\label{ensembles de mesure nulle}


Dans ce paragraphe, $\mathbf{D}$ désignera $\mathbf{D}_\!X$. Soit $\mathbf{D}^\ast$ la famille des  sous-ensembles $E\subset\G(\!X)$ pour lesquels il existe des ensembles mesurables $A$ et $B$ tels que $A\subset E\subset B$ et $\mu(B\backslash A)=0$.

\begin{lem}
Les ensembles $\mathbf{D}$ et $\mathbf{D}^\ast$ sont égaux.
\end{lem}

\begin{proof}
L'inclusion $\mathbf{D}\subset\mathbf{D}^\ast$ est claire. Soit $E\in\mathbf{D}^\ast$. Il s'agit de montrer que pour tout $\varepsilon>0$ on peut trouver une $\varepsilon$-approximation cylindrique $E_{M(\varepsilon)}(\varepsilon)$ de $E$.

Soit $\varepsilon>0$. Par définition, il existe des ensembles mesurables $A\in\mathbf{D}$ et $B\in\mathbf{D}$ tels que $A\subset E\subset B$ et $\mu(B\backslash A)=0$. Il existe un entier $n_0\in\mathbf{N}$, une $1/n_0$-approximation cylindrique $A_I$ de $A$ et une $1/n_0$-approximation cylindrique $B_J$ de $B$ telles que 

\begin{enumerate}
\item $1/n_0<\varepsilon$

\item  $\lVert\mu(A_i(1/n_0))\rVert<\varepsilon$ et $\lVert\mu(B_j(1/n_0))\rVert<\varepsilon$ pour tout $i\in I$ et tout $j\in J$.

\item $\lVert\mu(B_0(1/n_0)\backslash A_0(1/n_0))\rVert<\varepsilon$.

\end{enumerate}
La suite de cylindres définie par $E_0(\varepsilon):=B_0(1/n_0)$, par $B_0(1/n_0)\backslash A_0(1/n_0)$, par les $A_i(1/n_0)$ pour tout $i\in I$ et les $B_i(1/n_0)$ pour tout $i\in J$ est une $\varepsilon$-approximation de $E$.

\end{proof}

\begin{lem}
Soit $E\subset\G(\!X)$. Alors $E$ est mesurable de mesure nulle si et seulement si pour tout $\varepsilon>0$, il existe un ensemble au plus dénombrable $I$ et une suite de cylindres $(E_i(\varepsilon))_{i\in I}$ telle que $E\subset\cup_{i\in I} E_i(\varepsilon)$ et $\lVert\mu(E_i(\varepsilon))\rVert<\varepsilon$ pour tout $i\in I$.
\end{lem}

\begin{proof}
Cette condition est suffisante de manière évidente. Inversement, soit $E$ un ensemble mesurable de mesure nulle. Soit $\varepsilon>0$. Il s'agit de trouver une $\varepsilon$-approximation cylindrique $E_{I(\varepsilon)}$ de $E$ de partie principale vide.

Par construction de la mesure, il existe un entier $n_0\in\mathbf{N}$ et une $1/n_0$-approximation cylindrique $E'_{I}$ de $E$ tels que:

\begin{enumerate}
\item $1/n_0<\varepsilon$.

\item $\lVert\mu(E'_0(1/n_0))\rVert<\varepsilon$.

\item $\lVert\mu(E'_i(1/n_0))\rVert<\varepsilon$ pour tout $i\in I$.

\end{enumerate}

La suite définie par $E_0(\varepsilon):=\emptyset$ et par les $E'_{i}(1/n_0)$ pour tout $i\in I\cup\{0\}$ est une $\varepsilon$-approximation cylindrique de $E$ vérifiant la propriété annoncée.
\end{proof}

\begin{cor}
Soit $A$ un sous-ensemble de $\G(\!X)$, tel que $A\subset C$, avec $C$ un ensemble mesurable de mesure nulle de $\G(\!X)$. Alors $A$ est un ensemble mesurable de mesure nulle. 
\end{cor}

\begin{proof}
Par d\'efinition, il existe une suite de cylindres $(C_i(\varepsilon))_{i\in I}$ telle que $C\subset\cup_{i\in I} C_i(\varepsilon)$ et $\lVert\mu(C_i(\varepsilon))\rVert<\varepsilon$ pour tout $i\geq 1$. Posons $A_i(\varepsilon):=C_i(\varepsilon)$, pour tout $i\in I$. Comme $A\subset C$, $A\in\cup_{i\in I} A_i(\varepsilon)$.
\end{proof}


\section{Mesurabilité des images directes et inverses}


Dans cette section, nous allons étudier le problème de la stabilité de la propriété de mesurabilité par image directe et inverse sous un morphisme $h$.

Soit $h:\!Y\rightarrow\!X$ un $\mathbb{D}$-morphisme de schémas formels {admissibles} et réduits, tous deux de {pure dimension $d$}.


\subsection{Le jacobien d'un morphisme}
\label{Le jacobien d'un morphisme}


\subsubsection{La multiplicité d'un point le long d'un sous-schéma formel fermé}

Soit $\!Z\hookrightarrow\!X$ un sous-$R$-sch\'ema formel ferm\'e de $\!X$, d\'efini par un $\!O_\!X$-faisceau d'id\'eaux $\!I_\!Z$. Soit $x\in\G(\!X)$. Soit $\varphi:\Spf R_{k'}\rightarrow \!X$ le $\mathbb{D}$-morphisme qui correspond à $x$. On dispose du morphisme canonique de $\!O_\!X$-modules:
$$
\varphi^{*}\!I_\!Z\rightarrow\!O_{R_{k'}}.
$$
L'image de ce morphisme est un $\!O_{R_{k'}}$-faisceaux d'id\'eaux de $\!O_{R_{k'}}$. Deux cas se pr\'esentent alors: soit cet id\'eal est nul, soit il existe un entier naturel $n\in\mathbf{N}$ tel que $\varphi^{*}\!I_\!Z.\!O_{R_{k'}}=\pi^n.\!O_{R_{k'}}$. 

On appelle \textit{multiplicité de $x$ le long de $\!Z$} l'entier d\'efini par la valuation du g\'en\'erateur de l'id\'eal $\varphi^{*}\!I_\!Z.\!O_{R_{k'}}$ de $R_{k'}$ si cet id\'eal est non nul, l'infini sinon. On d\'efinit alors l'application $\textrm{mult}_\!Z(x):\G(\!X)\rightarrow\mathbf{N}\cup\{\infty\}$ comme l'application qui \`a $x\in\G(\!X)$ associe la multiplicité de $x$ le long de $\!Z$.

Soit $j:\!U\hookrightarrow\!X$ un sous-$R$-schéma formel ouvert, quasi-compact. Soit $x\in\G(\!X)$ tel que le $\mathbb{D}$-morphisme $\varphi$ correspondant rende commutatif le diagramme suivant
$$
\xymatrix{\Spf R_{k'}\ar[r]^{\varphi}\ar[d]_{\varphi_\!U}&\!X\\
\!U\ar[ur]_j\\}
$$
Les deux faisceaux d'idéaux de $\!O_{R_{k'}}$, $\varphi_\!U^\ast\!I_{\!Z\cap\!U}.\!O_{R_{k'}}$ et $\varphi^\ast\!I_{\!Z}.\!O_{R_{k'}}$, sont égaux (par commutativité du diagramme). En particulier, on en déduit la formule:
$$
\textrm{mult}_\!Z(x)=\textrm{mult}_{\!Z\cap\!U}(x).
$$ 

\begin{lem}
Soit $\!Z\hookrightarrow\!X$ un sous-$R$-sch\'ema formel ferm\'e de $\!X$. On a les deux assertions suivantes:
\begin{enumerate}

\item  Pour tout $n\in\mathbf{N}$, la fibre $\mathrm{mult}_\!Z^{-1}(n)$ est un ensemble cylindrique de $\G(\!X)$.

\item Pour la fibre de l'infini, on a l'\'egalit\'e:
$$
\mathrm{mult}_\!Z^{-1}(\infty)=\G(\!Z).
$$ 
\end{enumerate}
\end{lem}

\begin{proof}
On peut supposer que $\!X$ est un schéma formel affine. Si $\!X:=\Spf R\{x_1,\ldots,x_N\}/\!I$, notons $A:= R\{x_1,\ldots,x_N\}/\!I$. Le sous-schéma formel $\!Z$ est alors défini par la donnée d'un idéal $\!J$ de $A$. Dans ce contexte, l'idéal $\!I_\!Z.\!O_{R_{k'}}$ est simplement l'idéal de $R_{k'}$ engendré par les $g((\varphi_i)_{1\leq i\leq N})$ pour tout $g\in \!J$ et où le $N$-uplet de $R_{k'}^N$, $(\varphi_i)_{1\leq i\leq N}$, correspond au morphisme $\varphi$. \textbf{1}. Posons
$$
C_{\geq n}:=\pi_{\!Y,n}^{-1}(\G_n(\!Z)).
$$
On va montrer que
$$
C_{\geq n}=\{x\in\G(\!X)\mid\mathrm{mult}_\!Z(x)\geq n\}.
$$
Cette égalité découle alors directement des définitions. Comme $C_{\geq n}$ est un cylindre de $\G(\!X)$, la première assertion est une conséquence du fait que
$$
\{x\in\G(\!X)\mid \mathrm{mult}_\!Z(x)=n\}=C_{\geq n}\backslash C_{\geq n+1}.
$$
Le point \textbf{2}. se déduit de la description ci-dessus.
\end{proof}

\subsubsection{La d\'efinition}

Soit $y$ un point de $\G(\!Y)\backslash\G(\!Y_{\mathrm{sing}})$ et $\varphi:\Spf R_{k'}\rightarrow \!Y$ le morphisme correspondant, avec $k'$ une clôture parfaite de $\kappa(y)$. On dispose de la suite exacte canonique:
$$
h^{*}\Omega_{\!X/R}^1\rightarrow\Omega_{\!Y/R}^1\rightarrow\Omega_{\!Y/\!X}^1\rightarrow 0
$$ 

qui induit, en passant aux puissances ext\'erieures $d$-èmes, une nouvelle suite exacte:
$$
h^{*}\Omega_{\!X/R}^d\rightarrow\Omega_{\!Y/R}^d\rightarrow\Omega_{\!Y/\!X}^d\rightarrow 0.
$$
Par exactitude à droite du foncteur image inverse, on déduit une nouvelle suite exacte dans la catégorie des $\!O_{R_{k'}}$-Modules (équivalente à celle des $R_{k'}$-modules $\pi$-adiques, puisque $\Spf R_{k'}$ est affine):
$$
\varphi^{*}(h^{*}\Omega_{\!X/R}^d)\rightarrow\varphi^{*}\Omega_{\!Y/R}^d\rightarrow\varphi^{*}\Omega_{\!Y/\!X}^d\rightarrow 0.
$$
Le $R_{k'}$-morphisme $\varphi^{*}(h^{*}\Omega_{\!X/R}^d)\rightarrow\varphi^{*}\Omega_{\!Y/R}^d$ induit un morphisme:
$$
(\varphi^{*}(h^{*}\Omega_{\!X/R}^d))/(\mathrm{torsion})\rightarrow(\varphi^{*}\Omega_{\!Y/R}^d)/(\mathrm{torsion}).
$$
On définit alors l'\textit{ordre de l'idéal jacobien de $h$ en $y$}, noté $\ord_\pi(\mathrm{Jac})_h(y)$, de la manière suivante. Le choix de $\varphi$ entra\^\i ne que $L:=(h^{*}\Omega_{\!X/R}^d)/(\mathrm{torsion})$ est un $R_{k'}$-module libre de rang 1. On déduit du théorème de structure des modules de type fini sur les anneaux principaux, que l'image de $M:=(\varphi^{*}(h^{*}\Omega_{\!X/R}^d))/(\mathrm{torsion})$ dans $L$ est soit 0, soit $\pi^nL$, pour un certain $n\in\mathbf{N}$. On pose alors $\ord_\pi(\mathrm{Jac})_h(y)=\infty$ et $\ord_\pi(\mathrm{Jac})_h(y)=n$ respectivement.

\begin{defi}
On appelle ordre de l'idéal jacobien de $h$ l'application 
$$
\ord_\pi(\mathrm{Jac})_h:\G(\!Y)\backslash\G(\!Y_{\mathrm{sing}})\rightarrow\mathbf{N}\cup\{\infty\}
$$ 
qui à $y\in\G(\!Y)\backslash\G(\!Y_{\mathrm{sing}})$ associe l'entier $\ord_\pi(\mathrm{Jac})_h(y)$ défini ci-dessus.
\end{defi}

On appelle \textit{lieu sauvage de $h$} le sous-ensemble de $\G(\!Y)$ défini par
$$
\Sigma_h:=\G(\!Y_{\mathrm{sing}})\cup h^{-1}(\G(\!X_{\mathrm{sing}}))\cup\{y\in\G(\!Y)\mid \ord_\pi(\mathrm{Jac})_h(y)=\infty\}.
$$
On notera $\Pi_h$ l'ensemble $\{y\in\G(\!Y)\mid \ord_\pi(\mathrm{Jac})_h(y)=\infty\}$. Soit $B$ un sous-ensemble de $\G(\!Y)$. On dit que le morphisme $h$ est \textit{tempéré sur $B$} si $B\cap \Sigma_h$ est mesurable de mesure nulle. On dit que $h$ est \textit{tempéré} si $h$ est tempéré sur $\G(\!Y)$, \textit{i.e.} si le lieu sauvage de $h$ est de mesure nulle. 

Supposons que le $R$-schéma formel $\!Y$ est lisse. Soit $B$ un sous-ensemble de $\G(\!Y)$ contenu dans le complémentaire du lieu sauvage de $h$. L'ensemble $B$ est alors contenu dans l'ind-cylindre
$$
(\cup_{e\in\mathbf{N}}h^{-1}(\G^{(e)}(\!X))\cap(\cup_{e'\in\mathbf{N}}\{y\in\G(\!Y)\mid \ord_\pi(\mathrm{Jac})_h(y)=e'\}).
$$
Nous noterons $\Delta_{e,e'}$ l'ensemble $h^{-1}(\G^{(e)}(\!X))\cap\{y\in\G(\!Y)\mid \ord_\pi(\mathrm{Jac})_h(y)=e'\}$. En particulier, si $B$ est un cylindre, il est contenu dans la réunion d'un nombre fini de tels $\Delta_{e,e'}$.

\subsubsection{{Le cas lisse}.}

Supposons en outre que le $R$-sch\'ema formel $\!Y$ est lisse sur $R$. Dans ce cas, l'application $\ord_\pi(\mathrm{Jac})_h$ est définie (grâce à un procédé général) par un $\!O_\!Y$-faisceau d'idéaux, ${\Fitt}_d(\Omega_{\!Y/\!X}^1)$.

Soit $y\in\G(\!Y)$ et $\varphi\in\!Y(R_{k'})$ le $R$-morphisme qui lui correspond. L'image du morphisme naturel:
$$
M:=(\varphi^\ast(h^\ast\Omega^1_{\!X/R}))/(\mathrm{torsion})\rightarrow L:=(\varphi^\ast\Omega^1_{\!X/R})
$$ 

est un $\!O_{R_{k'}}$-Module, $M'$, libre de rang au plus $d$. La suite exacte canonique:
$$
0\rightarrow M'\rightarrow L\rightarrow \varphi^\ast\Omega^1_{\!Y/\!X}\rightarrow 0
$$
est une présentation de $\varphi^\ast\Omega^1_{\!Y/\!X}$, puisque $\!Y$ est lisse. De cette suite, on déduit le morphisme:
$$
\wedge^d(M')\otimes_{\!O_{R_{k'}}}\wedge^d L^\vee\rightarrow\!O_{R_{k'}}
$$
dont l'image est le $d$-ème idéal de Fitting de $\varphi^\ast\Omega^1_{\!Y/\!X}$.

\begin{lem}
Si $\!Y$ est lisse, les applications $\ord_\pi(\mathrm{Jac})_h$ et $\mathrm{mult}_{\!Z_h}$ sont \'egales sur $\G(\!Y)$, où $\!Z_h\hookrightarrow\!X$ est le sous-$R$-schéma formel de $\!X$ défini par le $d$-ème idéal de Fitting de $\Omega^1_{\!Y/\!X}$.
\end{lem}

\begin{proof}
Ceci découle directement des définitions des deux fonctions.
\end{proof}

\subsubsection{Une formule de composition.}

\begin{lem}
\label{formule de composition}
Soient $h:\!Y\rightarrow\!X$ et $g:\!Z\rightarrow\!Y$ deux $R$-morphismes de $R$-sch\'emas formels admissibles et de pure dimension $d$. Soit $z\in\G(\!Z)$ tel que
$$
z\not\in(\G(\!Z_{\mathrm{sing}})\cup g^{-1}(\G(\!Y_{\mathrm{sing}}))\cup (h\circ g)^{-1}(\G(\!X_{\mathrm{sing}}))).
$$
On a alors l'égalité suivante:
$$
\ord_\pi(\mathrm{Jac})_{h\circ g}(z)=\ord_\pi(\mathrm{Jac})_h\circ g(z)+\ord_\pi(\mathrm{Jac})_g(z).
$$
\end{lem}

\begin{proof}
Soit $\varphi\in\G(\!Z)(k')$ un $R$-morphisme correspondant à $z$. L'ordre de l'idéal jacobien de $h\circ g$ est donn\'e par le morphisme:
$$
(h\circ g)^{*}\Omega_{\!X/R}^d\rightarrow\Omega_{\!Z/R}^d.
$$ 
La fonction $\ord_\pi(\mathrm{Jac})_h$ est d\'ecrite par le morphisme: 
$$
h^{*}\Omega_{\!X/R}^d\rightarrow \Omega_{\!Y/R}^d.
$$ 
Enfin, l'ordre de l'idéal jacobien du morphisme $g$ est donn\'e par:
$$
g^{*}\Omega_{\!Y/R}^d\rightarrow \Omega_{\!Z/R}^d.
$$ 
Ces morphismes induisent deux $\!O_{R_{k'}}$-morphismes
$$
\varphi^{*}(h\circ g)^{*}\Omega_{\!X/R}^d\rightarrow\varphi^{*}\Omega_{\!Z/R}^d$$ 
et 
$$
\xymatrix{(g\circ\varphi)^{\ast}(h^{\ast}\Omega_{\!X/R}^d)\ar[r]& (g\circ\varphi)^{\ast}\Omega_{\!Y/R}^d\ar[r]^{=}&\varphi^{\ast}(g^{\ast}\Omega_{\!Y/R}^d)\ar[r]& \varphi^{\ast}\Omega_{\!Z/R}^d}
$$
qui sont égaux. Ceci entraine en particulier que l'ordre de l'idéal jacobien de $h\circ g$ en $z$ est infini si et seulement si l'ordre de l'idéal jacobien de $h$ est infini en $g(z)$ ou si l'ordre du jacobien de $g$ est infini en $z$.

Supposons désormais que l'ordre de l'idéal jacobien de $h\circ g$ en $z$ est fini. Posons $L:=(g\circ\varphi)^{\ast}\Omega_{\!Y/R}^d/(\mathrm{torsion})$. L'image de $M:=(h\circ g \circ\varphi)^{\ast}\Omega_{\!X/R}^d/(\mathrm{torsion}$) dans $L$ est un $\!O_{R_{k'}}$-Module $M'$ libre de type fini. Il existe $n\in\mathbf{N}$ tel que $M'=\pi^n L$. De même l'image $L'$ de $L$ dans $K:=\varphi^{*}\Omega_{\!Z/R}^d/(\mathrm{torsion})$ est canoniquement isomorphe à $L'=\pi^m K$. L'égalité des deux morphismes précédents prouve la formule:
$$
\ord_\pi(\mathrm{Jac})_{h\circ g}(y)=m+n
$$
qui est l'égalité voulue.

\end{proof}


\subsection{La mesurabilit\'e de l'image directe}
\label{mesurabilite et image directe}


Soit $B$ un ensemble mesurable de $\G(\!Y)$. Dans ce paragraphe, nous étudions le problème de la mesurabilité dans $\G(\!X)$ de l'image par $h$ de l'ensemble $B$.

\subsubsection{Le cas cylindrique}

En général, l'image par $h$ d'un cylindre de $\G(\!Y)$ n'est pas un cylindre de $\G(\!X)$ (même si $\!Y$ est supposé lisse).

\begin{lem}
\label{lemme de mesurabilite et image directe}
Supposons que $\!Y$ est lisse. Si $B$ est un cylindre de $\G(\!Y)$, alors on a les propri\'et\'es suivantes:

\begin{enumerate}

\item il existe une suite de cylindres $(A_i)_{i\in\mathbf{N}}$ de $\G(\!X)$ telle que $h(B)\subset \cup_{i\in\mathbf{N}} A_i$ et $\lVert\mu_\!X(A_i)\rVert\leq\lVert\mu_\!Y(B)\rVert$, pout tout $i\in\mathbf{N}$.

\item Supposons que $h(B)\cap\G(\!X_{\mathrm{sing}})=\emptyset$, et que $B\cap\Pi_h=\emptyset$. Alors $h(B)$ est un cylindre.
\end{enumerate}
\end{lem}

\begin{proof}
\textbf{1}. Supposons d'abord que $\lVert\mu_\!Y(B)\rVert=0$. Comme $B$ est un cylindre et que $\!Y$ est lisse, ceci entra\^\i ne que $B=\emptyset$ et que $h(B)=\emptyset$. Supposons désormais que $\lVert\mu_\!Y(B)\rVert\not=0$. Comme $\G(\!X_{\mathrm{sing}})$ est contenu dans un cylindre de volume aussi petit que l'on veut, on peut supposer que $h(B)\subset\G(\!X)\backslash\G(\!X_{\mathrm{sing}})$. Soit $e_0\in\mathbf{N}$ un entier. On a l'égalité ensembliste
$$
h(B)=\cup_{e\geq e_0}(\G^{(e)}(\!X)\cap h(B)).
$$
Posons $\G^{(e)}(\!X)\cap h(B)=:h^{(e)}(B)$ pour tout $e\geq e_0$. Soit $e\geq e_0$. Notons $k_B$ le rang du cylindre $B$. Pour $n\geq \max(k_B,e_0,c_\!Xe)$, l'ensemble
$$
A^{(e)}_n:=\pi_{n,\!X}^{-1}(\pi_{n,\!X}(h^{(e)}(B)))
$$
est un cylindre de $\G(\!X)$ par le théorème de Chevalley, de rang $n$ et contenant $h^{(e)}(B)$. En outre, $A^{(e)}_n$ est contenu dans $\G^{(e)}(\!X)$. Soit $n\geq \max(k_B,e_0,c_\!Xe)$, comme $A^{(e)}_n$ est stable au rang $n$, 
$$
\mu_\!X(A^{(e)})=[\pi_{n,\!X}(h^{(e)}(B))]\mathbf{L}^{-(n+1)d} \ \ \textrm{et} \ \ \lVert\mu_{\!X}(A^{(e)})\rVert\leq 2^{-(n+1)d+\dim\pi_{n,\!X}(B)}.
$$
Or  $\lVert\mu_\!Y(B)\rVert=2^{-(n+1)d+\dim\pi_{n,\!X}(B)}$.

\textbf{2}. Par le lemme \ref{quasi-compacite de la topologie constructible}, il existe un entier $e'$ tel que $B$ est contenu dans $h^{-1}(\G^{(e')}(\!X))$ et l'on peut supposer que $\ord_\pi(\mathrm{Jac})_h$ est borné sur $B$. L'assertion découle alors du lemme clé \ref{5}.
\end{proof}

\begin{cor}
\label{bidouille}
Supposons que $\!Y$ est lisse. Soient $B\subset\G(\!Y)$ et $A\subset\G(\!X)$ des ensembles fortement mesurables tels que $B\cap\Sigma_h=\emptyset$ et $h(B)=A$. Alors, pour tout $\varepsilon>0$, il existe une $\varepsilon$-approximation cylindrique $B_{J(\varepsilon)}$ de $B$ et une $\varepsilon$-approximation cylindrique $A_{I(\varepsilon)}$ de $A$ telles que $h(B_0(\varepsilon))=A_0(\varepsilon)$.
\end{cor}

\begin{proof}
Soit $B_J$ une $\varepsilon$-approximation cylindrique de $B$. Par hypothèse, $B=\cup_{e,e'\in\mathbf{N}} (B\cap \Delta_{e,e'})$ et $B_0(\varepsilon)=\cup_{e,e'\in\mathbf{N}} (B_0(\varepsilon)\cap \Delta_{e,e'})$. Posons $B_{e,e'}:= B\cap \Delta_{e,e'}$ et $B_i^{(e,e')}(\varepsilon):= B_i(\varepsilon)\cap \Delta_{e,e'}$, pour tout $i\in J\cup\{0\}$. Par le lemme \ref{5}, l'ensemble $h(B_0^{(e,e')}(\varepsilon))$ est un cylindre de $\G(\!X)$. En outre, le lemme \ref{quasi-compacite de la topologie constructible} assure qu'il existe $m\in\mathbf{N}$ tel que
$$
B_0(\varepsilon)=\cup_{1\leq i,j\leq m}B_0^{(e_i,e'_j)}(\varepsilon).
$$
L'égalité $B\triangle B_0(\varepsilon)\subset \cup_{i\geq 1}A_i$ entra\^\i ne que, pour tout $e$ et $e'\in\mathbf{N}$, $B_{e,e'}\triangle B_0^{(e,e')}(\varepsilon)\subset \cup_{i\geq 1}B_i^{(e,e')}$ et donc que $h(B_{e,e'})\triangle h(B_0^{(e,e')}(\varepsilon))\subset \cup_{i\geq 1}h(B_i^{(e,e')})$. Comme 
$$
\cup_{e,e'\in\mathbf{N}} h(B_0^{(e,e')}(\varepsilon))=\cup_{1\leq i,j\leq m}h(B_0^{(e_i,e'_j)}(\varepsilon)),
$$ 
l'ensemble $h(B_0(\varepsilon))$ est un cylindre de $\G(\!X)$. Le lemme \ref{lemme de mesurabilite et image directe}/1 assure la validité de l'assertion annoncée.
\end{proof}

\subsubsection{Le cas général}

\begin{thm}
\label{theoreme de image directe}
Supposons que $\!Y$ est lisse et que $h$ est tempéré. Si $B\subset\G(\!Y)$ est mesurable (resp. fortement mesurable), alors $h(B)\subset \G(\!X)$ est mesurable (resp. fortement mesurable).
\end{thm}

\begin{proof}
Soit $\varepsilon>0$. Posons $A:=h(B)$. On a
$$
B=(B\cap h^{-1}(\G(\!X_{\mathrm{sing}})))\sqcup(B\backslash h^{-1}(\G(\!X_{\mathrm{sing}}))).
$$
Comme $h$ est tempéré, on peut supposer que $B\cap\Sigma_h=\emptyset$. Soit $B_{I(\varepsilon)}$ (\textit{resp.} $\Sigma_{J(\varepsilon)}$) une $\varepsilon$-approximation cylindrique de $B$ (\textit{resp.} $\Sigma_h$). Comme $B\cap\Sigma_h=\emptyset$, on peut supposer que $B_0(\varepsilon)\cap \Sigma_0(\varepsilon)=\emptyset$. Comme l'ensemble $\Sigma_h$ est un pro-constructible de $\G(\!Y)$, la quasi-compacité de la topologie constructible entraine que
$$
\Sigma_h\subset(\Sigma_0(\varepsilon))\cup(\cup_{1\leq i\leq m}\Sigma_i(\varepsilon)),
$$
avec $m\in\mathbf{N}$. Comme la suite de cylindres de $\G(\!Y)$, définie par 
$$
B'_0(\varepsilon):=B_0(\varepsilon)\backslash (\cup_{1\leq i\leq m}\Sigma_i(\varepsilon)),
$$
par les $= \Sigma_i(\varepsilon)$ pour $1\leq i\leq m$ et les $B_i(\varepsilon)$ pour $i\in I(\varepsilon)$, est une $\varepsilon$-approximation cylindrique de $B$, on peut supposer que $B_{I(\varepsilon)}$ vérifie la condition $B_0(\varepsilon)\cap\Sigma_h=\emptyset$. Or
$$
h(B)\triangle h(B_0(\varepsilon))\subset \cup_{i\in I(\varepsilon)}h(B_i(\varepsilon)).
$$
et, si $B$ est fortement mesurable, $B_0(\varepsilon)\subset B$ et $h(B_0(\varepsilon))\subset A$. Le théorème est alors une conséquence directe du lemme \ref{lemme de mesurabilite et image directe}.
\end{proof}


\subsection{La mesurabilit\'e de l'image inverse}
\label{mesurabilite et image inverse}


Soit $A$ une partie mesurable de $\G(\!X)$. Dans ce paragraphe, nous étudions le problème de la mesurabilité de l'image inverse $h^{-1}(A)$ par le morphisme $h$ de l'ensemble $A$.

\subsubsection{Le cas cylindrique.}

Supposons que $A\in\mathbf{D_\!X}$ soit un sous-ensemble cylindrique de $\G(\!X)$.

\begin{lem}
Soit $A$ un cylindre de $\G(\!X)$ de rang $n$. Alors $h^{-1}(A)$ est un cylindre de $\!Y$ de rang $n$.
\end{lem}

\begin{proof}
Cela d\'ecoule directement de la d\'efinition d'un cylindre et de la commutativit\'e du diagramme:
$$
\xymatrix{\G(\!Y)\ar[r]^{h}\ar[d]_{\pi_{n,\!Y}}&\G(\!X)\ar[d]^{\pi_{n,\!X}}\\
\G_n(\!Y)\ar[r]_{h_n}&\G_n(\!X)}
$$ 
où $h_n$ désigne l'image par le foncteur $\G_n$ du $R_n$-morphisme de schémas $h_n:Y_n\rightarrow X_n$, induit par $h$. 
\end{proof}

\subsubsection{Le cas injectif}
\label{hiv pour tcvg}

Supposons, dans ce paragraphe, que le $k$-morphisme de schémas $h:\G(\!Y)\rightarrow\G(\!X)$ est injectif et tempéré et que $\!Y$ est lisse.

\begin{exe}
Soit $\!X$ un $R$-schéma formel admissible, dont la fibre générique $X_K$ est lisse. La théorie des modèles de Néron faibles (cf \cite{blr}, \cite{bs} et \cite{ls} \S 2.6) assure que l'on peut trouver un $R$-schéma formel admissible $\!U$, lisse sur $R$, (cf \ref{neron=sttf}) et un $\mathbb{D}$-morphisme de schémas formels $h:\!U\rightarrow\!X$ tels que $h:\G(\!U)\rightarrow\G(\!X)$ soit tempéré et fortement injectif (\textit{cf} \S\ref{fi}).
\end{exe}

\begin{thm}
\label{hiv pour tcvg 1}
Sous les hypothèses ci-dessus, si $A$ est mesurable (resp. fortement mesurable) dans $\G(\!X)$, alors l'ensemble $B:=h^{-1}(A)$ est mesurable (resp. fortement mesurable) dans $\G(\!Y)$.
\end{thm}

\begin{proof}
Soit $A_{I(\varepsilon)}$ une $\varepsilon$-approximation de $A$ telle que $A_0(\varepsilon)\subset A$. Comme 
$$
A=(A\cap\G(\!X_{\mathrm{sing}}))\sqcup(A\backslash\G(\!X_{\mathrm{sing}})),
$$
et comme $h$ est tempéré, on peut supposer que $A\cap\G(\!X_{\mathrm{sing}})=\emptyset$ et que $B\cap\Pi_h=\emptyset$. La définition de $A_{I(\varepsilon)}$ entraine que
$$
B\triangle h^{-1}(A_0(\varepsilon))\subset \cup_{i\in I(\varepsilon)}h^{-1}(A_i(\varepsilon)).
$$
Soit $1\geq\delta>0$. Par hypothèse, le lieu sauvage de $h$ est de mesure nulle. Il existe donc un cylindre $D(\delta)\supset \Sigma_h$ tel que:
\begin{enumerate}

\item $\lVert D(\delta)\rVert<\delta$.

\item $\Sigma_h=\cap_\delta D(\delta)$.
\end{enumerate}
On a donc
$$
(B\backslash D(\delta))\triangle (h^{-1}(A_0(\varepsilon))\backslash D(\delta))\subset \cup_{i\in I(\varepsilon)}(h^{-1}(A_i(\varepsilon))\backslash D(\delta)).
$$
Par définition de $\Sigma_h$ et grâce au lemme \ref{quasi-compacite de la topologie constructible}, il existe un entier $m:=m(\delta)\geq 1$ tel que $\G(\!Y)\backslash D(\delta)\subset \cup_{1\leq i, j\leq m}\Delta_{e_i,e_j}$. La suite de cylindres de $\G(\!Y)$, définie par $(h^{-1}(A_i(\varepsilon))\backslash D(\delta))$, est une $(r(\delta)\times\varepsilon)$-approximation cylindrique de $B$, où $r(\delta):=\max(2^{e_j})$ est donné par le lemme clé \ref{5}. Pour tout $\delta>0$, $B\backslash D(\delta)$ est un ensemble mesurable de $\G(\!Y)$.

Il nous faut en déduire que $B$ lui-même est mesurable. On peut supposer que pour $0<\delta'<\delta\leq 1$, $D(\delta')\subset D(\delta)$. Le résultat précédent entra\^\i ne alors que $B\cap (D(\delta)\backslash D(\delta'))$ est mesurable. En outre, comme par hypothèse $\mu_\!Y(\Sigma_h)=0$, la famille $\mu_\!Y(D(\delta))$ tend vers 0. L'égalité
$$
(B\backslash D(1))\sqcup(\sqcup_{0<\delta'<\delta\leq 1}(B\cap D(\delta)\backslash D(\delta'))=B\backslash \Sigma_h=B
$$
permet de conclure (\textit{cf} proposition \ref{mesurable comme limite}).

\end{proof}


\section{Intégrale motivique et théorèmes de changement de variables}
\label{fonctions integrables}



\subsection{La définition de l'intégrale}


On supposera que tous les $R$-schémas formels sont admissibles et réduits.

\begin{defi}
\label{definition de fonction mesurable}
On dit qu'une application $\psi:\G(\!X)\rightarrow\widehat{\!M}$ est mesurable si les fibres de $\psi$ sont mesurables dans $\G(\!X)$ et si la somme $\sum_{a\in\widehat{\!M}}\mu_\!X(\psi^{-1}(a))\psi(a)$ est convergente, \textit{i.e.} si l'ensemble des $\mu_\!X(\psi^{-1}(a))$ non nuls pour ${a\in\widehat{\!M}}$ est au plus d\'enombrable et si $\sum\mu_\!X(\psi^{-1}(a_i))a_i$ converge dans $\widehat{\!M}$.
\end{defi}

\begin{defi}
Si $\psi$ est une fonction int\'egrable sur $\G(\!X)$, on d\'efinit l'int\'egrale motivique de $\varphi$ par: 
$$
\int_{\G(\!X)}\psi d\mu_\!X:=\sum_{i\in \mathbf{N}}\mu_\!X(\psi^{-1}(a_i))a_i.
$$
Si $A$ est un ensemble mesurable de $\G(\!X)$ et $\psi:A\rightarrow\widehat{\!M}$ une application v\'erifiant les propri\'et\'es de la d\'efinition \ref{definition de fonction mesurable}, on d\'efinit de m\^eme l'int\'egrale de $\psi$ sur $A$ par:
$$
\int_{A}\psi d\mu_\!X:=\sum_{i\in \mathbf{N}}\mu_\!X(\psi^{-1}(a_i))a_i.
$$
\end{defi}

\begin{defi}
Pour tout ensemble mesurable $A$ et toute application $\alpha:A\rightarrow\mathbf{Z}\cup\lbrace\infty\rbrace$, on dit que $\alpha$ est exponentiellement int\'egrable si $\mathbf{L}^{-\alpha}$ est int\'egrable, \textit{i.e.} si les fibres finies de $\alpha$ sont mesurables et si la somme $\Sigma_{n\in\mathbf{Z}}\mu_\!X(\alpha^{-1}(n))\mathbf{L}^{-n}$ converge dans $\widehat{\!M}$, ce qui équivaut au fait que la suite $(\lVert\mu_\!X(\alpha^{-1}(n))\rVert/2^n)_{n\in\mathbf{N}}$ converge vers 0 dans $\widehat{\!M}$.
\end{defi}

\begin{lem}
Soient $A$ et $B$ deux parties mesurables de $\G(\!X)$ telles que $A\subset B$. Si $\psi:B\rightarrow\widehat{\!M}$ est une application int\'egrable, alors sa restriction \`a $A$, $\psi_{\mid A}:A\rightarrow \widehat{\!M}$ est int\'egrable sur $A$.
\end{lem}

\begin{proof}
On a la relation $\psi^{-1}(a)\cap A\subset \psi^{-1}(a)$ pour tout $a\in\widehat{\!M}$. En particulier, ceci entra\^\i ne que l'ensemble des $((\psi_{\mid A})^{-1}(a))_{a\in\widehat{\!M}}$ non nuls est au plus d\'enombrable. Comme $\widehat{\!M}$ est complet et ultram\'etrique, la s\'erie $\sum_{i\in\mathbf{N}}(\psi^{-1}(a_i)\cap A)a_i$ converge si et seulement si le terme g\'en\'eral tend vers 0. Soit $\varepsilon>0$, il existe $N\in\mathbf{N}$ tel que pour tout $n\geq N$,  $\lVert\mu(\psi^{-1}(a_n))\rVert<\varepsilon$. L'inclusion du d\'ebut entraine pour ces m\^emes indices, que $\lVert\mu((\psi_{\mid A})^{-1}(a_n))\rVert<\varepsilon$. D'o\`u le r\'esultat.

\end{proof}

\begin{lem}
Soient $A$ et $B$ deux parties mesurables. Soit $\psi$ une application int\'egrable sur $A\cup B$, alors on a la relation:
$$
\int_{A\cup B} \psi d\mu_\!X=\int_A \psi d\mu_\!X +\int_B \psi d\mu_\!X -\int_{A\cap B} \psi d\mu_\!X.
$$
\end{lem}

\begin{proof}
Par d\'efinition, on a:
$$
\int_{A\cup B} \psi d\mu_\!X=\sum_{i\in \mathbf{N}}\mu_\!X(\psi^{-1}(a_i))a_i.
$$ 
Le r\'esultat vient alors du fait que $\mu_\!X(\psi^{-1}(a_i))=\mu_\!X(\psi^{-1}(a_i)\cap A) +\mu_\!X(\psi^{-1}(a_i)\cap B)-\mu_\!X(\psi^{-1}(a_i)\cap (A\cap B))$. Autrement dit, $\mu_\!X(\psi^{-1}(a_i))=\mu_\!X((\psi_{\mid A})^{-1}(a_i)) +\mu_\!X((\psi_{\mid B})^{-1}(a_i))-\mu_\!X((\psi_{\mid A\cap B})^{-1}(a_i))$. Le lemme pr\'ec\'edent conclut.
\end{proof}

\begin{lem}
Soit $A$ un ensemble mesurable de $\G(\!X)$. Soit $B$ un ensemble de mesure nulle. Soit $\psi$ une application int\'egrable sur $A$. On a alors la relation:
$$
\int_A\psi d\mu=\int_{A\backslash B}\psi d\mu.
$$
\end{lem}

\begin{proof}
La preuve ressemble \`a la pr\'ec\'edente. Il suffit de remarquer que $\mu(\psi^{-1}(a_i))=\mu(\psi^{-1}(a_i)\cap (A\cap B)) +\mu(\psi^{-1}(a_i)\cap (A\backslash B)$. Or $\psi^{-1}(a_i)\cap (A\cap B)$ est contenu dans $B$, donc de mesure nulle. Donc $\mu(\psi^{-1}(a_i))=\mu(\psi^{-1}(a_i)\cap (A\backslash B))$.
\end{proof}


\subsection{Les th\'eor\`emes de changement de variables}


Soit $h:\!Y\rightarrow\!X$ un $R$-morphisme de schémas formels admissibles et réduits, tous deux de pure dimension $d$.

\begin{lem}
\label{4}
Supposons que $\!Y$ est lisse.  Alors pour tout $n\geq \mathrm{max}(e,c_\!Xe')$, pour tout $z\in \Delta_{e,e'}$ et tout $x\in\G(\!X)$, tel que $\pi_{n,\!X}(h(z))=\pi_{n,\!X}(x)$, il existe $y\in\G(\!Y)$ tel que $h(y)=x$ et $\pi_{n-e,\!Y}(z)=\pi_{n-e,\!Y}(y)$.
\end{lem}

\begin{proof}
Soit $z\in\G(\!Y)$ tel que $\varphi_z\in \Delta_{e,e'}(k')$  et $x\in\G(\!X)$. Quitte à étendre les scalaires, on peut supposer que $\varphi_z:\mathbb{D}\rightarrow \!Y$ et qu'à $x$ correspond un $R$-morphisme $\varphi_x:\mathbb{D}\rightarrow\!X$. La question étant locale en $\!Y$, on peut supposer que $\!Y$ est affine et comme $\!Y$ est supposé lisse, on peut supposer qu'il existe un $R$-morphisme étale $\!Y\rightarrow\mathbb{B}^d_R$. On peut supposer que $\!Y=\mathbb{B}^d_R$. En effet, le lemme \ref{morphisme etale} assure qu'au-dessus de tout point la fibre du $k$-morphisme de schémas $\G(\!Y)\rightarrow (Y_0)_{\textrm{red}}$ est isomorphe à $\G(\mathbb{B}^d_R)$. Par ailleurs, la question est locale en $\!X$, on peut donc supposer que $\!X$ est affine.

\textit{Premier cas}: supposons que $\!X=\mathbb{B}^d_R$. Il nous suffit de prouver que pour tout $v\in R^d$, il existe $u\in R^d$ vérifiant l'équation
$$
h(\varphi_z+\pi^{n+1-e}u)=h(\varphi_z)+\pi^{n+1}v.\leqno (1)
$$
Ce système est alors équivalent à
$$
\pi^{-e}\!J_h(\varphi_z).u+\pi(\ldots)=v \leqno (2)
$$
où $\!J_h$ désigne la matrice jacobienne de $h$. Soit $M$ la matrice adjointe de $\!J_h$. On est ramené à résoudre le système
$$
\pi^{-e}M(\varphi_z)\!J_h(\varphi_z).u+\pi(\ldots)=M(\varphi_z)v. \leqno (3)
$$
Par définition de la fonction $\ord_\pi(\mathrm{Jac})_h$, comme $\!Y$ et $\!X$ sont lisses de même dimension, le déterminant de $\!J_h(\varphi_z)$ est de valuation exactement $e$. Par conséquent, la réduction modulo $\pi$ du système $(3)$ précédent est linéaire sur le corps $k$, de déterminant non nul. Il admet donc une solution par le lemme de Hensel (\textit{cf} \cite{bourb} III \S 4-3).

\textit{Deuxième cas}: $\!X=\Spf R\{x_1,\ldots,x_N\}/\!I$. Notons $\tilde{h}$ le $R$-morphisme composé
$$
\xymatrix{\mathbb{B}^d_R\ar[r]_h\ar@/^2pc/[rr]^{\tilde{h}}&\!X\ar@{^{(}->}[r]_{j}&\mathbb{B}^N_R}
$$
où $j$ est le morphisme d'immersion. 

La question étant locale en $\!X$, on peut supposer, comme dans la preuve du lemme \ref{2}, que $\!X$ est d'intersection complète, \textit{i.e.} qu'il existe $N-d$ éléments $(f_i)_{1\leq i\leq N-d}$ de $\!I$ qui l'engendrent. Soit $\Delta$ la matrice $(\partial f_j/\partial x_i)_{1\leq i\leq N, 1\leq j\leq N-d}$. Quitte à réordonner les $x_i$, on peut supposer qu'il existe un entier $e''\leq c_\!Xe'$ tel que la valuation du mineur $(N-d)\times(N-d)$, calculé à partir des $N-d$ premières colonnes et évalué en $\tilde{h}(\varphi_z)$, est exactement $e''$ et les valuations des autres mineurs $(N-d)\times(N-d)$ sont supérieures ou égales à $e''$.

Pour prouver l'assertion dans ce cas, il nous suffit de montrer que pour tout $v\in R^N$ tel que $f_j(h(\varphi_z) + \pi^{n+1}v)=0$, pour tout $1\leq j\leq N-d$, il existe $u\in R^d$ tel que
$$
h(\varphi_z+\pi^{n+1-e}u)=h(\varphi_z)+\pi^{n+1}v. \leqno{(1')}
$$
Soit $v\in R^N$. On peut se ramener au cas où $\!X=\mathbb{B}_R^d$. En effet, le $\mathbb{D}$-morphisme
\bigskip
$$
\xymatrix{\mathbb{B}_R^d\ar[r]^{\tilde{h}}\ar@/^2pc/[rr]&\mathbb{B}_R^N\ar[r]&\mathbb{B}^{N-d}_R}
$$
défini par
$$
\xymatrix{f_i(h_1(x_1,\ldots,x_d),\ldots,h_N(x_1,\ldots,x_d))&f_i(z_1,\ldots,z_N)\ar@{|->}[l]&y_i\ar@{|->}[l]}
$$
est constant d'image 0 (par définition de $\!X$). En particulier, ceci entra\^\i ne que le produit des matrices jacobiennes $\Delta(\tilde{h}(\varphi_z)).\!J_h(\varphi_z)$ est nul. Ceci implique que les colonnes de la matrice $\!J_h(\varphi_z)$ appartiennent au noyau de la matrice $\Delta(\tilde{h}(\varphi_z))$, \textit{i.e.} les coordonnées de chacune de ces colonnes sont solutions du système
$$
\pi^{-e''}M'(\tilde{h}(\varphi_z)).\Delta(\tilde{h}(\varphi_z))w=0 \leqno (2')
$$
où $M'$ est la matrice $N-d\times N-d$ définie comme la matrice adjointe de la sous-matrice $N-d\times N-d$ de $\Delta$ formée des $N-d$ premières colonnes et $w\in R^N$. Par hypothèse, les équations définissant le noyau de $\Delta(\tilde{h}(\varphi_z))$ expriment donc les $N-d$ premières coordonnées de $w$ en fonction des $d$ dernières. Le lemme de Hensel (\textit{cf} \cite{bourb} III \S 4-3) assure que $w\in R^N$ est solution du système $(2')$ si et seulement si le vecteur $\overline{w}\in k^N$, induit par $w\in R^N$ par réduction modulo $\pi$, est solution du système $(2')$ réduit modulo $\pi$. En particulier, comme pour tout $1\leq j\leq N-d$, $f_j(h(\varphi_z) + \pi^{n+1}v)=0$, $v$ est solution du système $(2')$. Le système $(1')$ étant équivalent à
$$
\pi^{-e}\!J_h(\varphi_z).u+\pi(\ldots)=v, \leqno (3')
$$
il suffit de résoudre
$$
\pi^{-e}\widetilde{\!J_h}(\varphi_z).\left(\begin{array}{c}
u_{N-d+1} \\ 
u_{N-d+2} \\ 
\vdots\\ 
u_N
\end{array}\right) +\pi(\ldots)=\left(\begin{array}{c}
v_{N-d+1} \\ 
v_{N-d+2} \\ 
\vdots\\ 
v_N
\end{array}\right)
$$
où $\widetilde{\!J_h}$ est la sous-matrice de ${\!J_h}$ formée des $d$ dernières lignes. Soient $x_{N-d+1},\ldots,x_{N}$  les $d$ dernières coordonnées parmi $x_1,\ldots,x_N$. Ces coordonnées définissent un $\mathbb{D}$-morphisme de schémas formels
$$
p:\!X\rightarrow\mathbb{B}_R^d.
$$
Par construction, la matrice jacobienne du morphisme $p\circ h$ n'est autre que la matrice $\widetilde{\!J_h}$. En outre, les hypothèses assurent que l'image de la base $(h\circ\varphi_z)^\ast dx_{i_{N-d+1}},\ldots,(h\circ\varphi_z)^\ast dx_{i_{N}}$ par le morphisme
$$
(h\circ\varphi_z)^\ast p^\ast\Omega^1_{\mathbb{B}_R^d/\mathbb{D}}\rightarrow (h\circ\varphi_z)^\ast\Omega^1_{\!X/\mathbb{D}}/(\mathrm{torsion})
$$
définit une base de $(h\circ\varphi_z)^\ast\Omega^1_{\!X/\mathbb{D}}/(\textrm{torsion})$. Autrement dit, $\ord_\pi(\mathrm{Jac})_p(h(\varphi_z))=0$. La formule de composition \ref{formule de composition} implique que $\ord_\pi(\mathrm{Jac})_{p\circ h}(\varphi_z)=e$. On est donc ramené au premier cas en considérant le $\mathbb{D}$-morphisme $p\circ h:\mathbb{B}_R^d\rightarrow \mathbb{B}_R^d$.
\end{proof}

\begin{lem}
\label{5}
Supposons que $\!Y$ est lisse. Soient $B\subset\G(\!Y)$ un cylindre et $A=h(B)$. Supposons que $\ord_{\pi}(\mathrm{Jac})_h(\varphi)$ est constant de valeur $e<\infty$ pour tout $\varphi\in B$, et $A\subset\G^{(e')}(\!X)$, avec $e'\in\mathbf{N}$. Alors $A$ est un cylindre. De plus, si la restriction de $h$ \`a $B$ est injective, alors, pour $n\in\mathbf{N}$ suffisamment grand, on a: 

\begin{enumerate}

\item Si $\varphi$ et $\varphi'$ appartiennent \`a $B$ et $\pi_{n,\!X}(h(\varphi))=\pi_{n,\!X}(h(\varphi'))$, alors $\pi_{n-e,\!Y}(\varphi)=\pi_{n-e,\!Y}(\varphi')$.

\item Le morphisme $h_n:\pi_{n,\!Y}(B)\rightarrow\pi_{n,\!X}(A)$ induit par $h$ est une fibration triviale par morceaux de fibre $\mathbf{A}_k^e$.
\end{enumerate}
\end{lem}

\begin{proof}
Prouvons tout d'abord que $A:=h(B)$ est un cylindre. Par hypothèse, $B$ est contenu dans $\Delta_{e,e'}$. Soit $n\in\mathbf{N}$ tel que $B$ soit un cylindre de rang $n-e$ (en particulier, $n\geq 2e$ et $n\geq 2e+c_\!Xe'$). Soit $x\in\pi_{n,\!X}^{-1}(\pi_{n,\!X}(A))$. Il existe donc $b\in B$ tel que
$$
\pi_{n,\!X}(x)=\pi_{n,\!X}(h(b)).
$$
Le lemme \ref{4} assure alors qu'il existe $y\in\G(\!Y)$ tel que
\begin{enumerate}

\item[(\textit{a})] $h(y)=x$.

\item[(\textit{b})] $\pi_{n-e,\!Y}(b)=\pi_{n-e,\!Y}(y)$.
\end{enumerate}
Comme $B$ est un cylindre de rang $n-e$, la condition 2 entraine que $y\in B$. La condition 1 assure que $x\in A$. On a donc montré que $A=\pi_{n,\!X}^{-1}(\pi_{n,\!X}(A))$. Par ailleurs, la relation $\pi_{n,\!X}(A)=h_n(\pi_{n,\!Y}(B))$, donc est constructible d'après un théorème déjà cité (\textit{cf} \cite{ega} \S 7). Ceci implique que l'ensemble $A$ est un cylindre. Prouvons l'assertion \ref{5}/1. Soit $n\in\mathbf{N}$ tel que $B$ soit un cylindre de rang $n-e$. Soient $\varphi$ et $\varphi'$ deux points de $B\subset\Delta_{e,e'}$ tels que $\pi_{n,\!X}(h(\varphi))=\pi_{n,\!X}(h(\varphi'))$. Le lemme \ref{4}, appliqué à $\varphi'$ et $\varphi$, assure l'existence de $y\in B$ tel que $y=\varphi'$ (par injectivité de $h$) et $\pi_{n-e,\!Y}(\varphi)=\pi_{n-e,\!Y}(y)=\pi_{n-e,\!Y}(y)$. Passons à l'assertion \ref{5}/2. Il nous suffit de démontrer le résultat pour le cylindre $\Delta_{e,e'}$. La question étant locale en $\!X$ et $\!Y$, on peut supposer que $\!X$ et $\!Y$ sont affines. Supposons que $\!X$ est défini dans $\mathbb{B}_R^N$ par l'idéal $\!I$, engendré par $m$ éléments $(f_i)_{1\leq i\leq m}$. Comme le cylindre $h(\Delta_{e,e'})$ est contenu dans $\G^{(e')}(\!X)$, on peut supposer (\textit{cf} lemme \ref{2}) qu'il est recouvert par un nombre fini de cylindres $C$, vérifiant la propriété
$$
\G(\!X)\cap C=\G(R\{X_1,\ldots,X_N\}/(f_{i_1},\ldots,f_{i_{N-d}}))\cap C
$$
On peut donc supposer que $\!X\hookrightarrow\mathbb{B}_R^N$ est d'intersection complète. Soit $\Delta$ la matrice jacobienne $(\partial f_j/\partial X_i)_{1\leq i\leq N,1\leq j\leq N-d}$. Quitte à restreindre $C$ et éventuellement à réordonner les $X_i$, on peut supposer que, pour tout $x\in C$, le mineur $N-d\times N-d$ de $\Delta(x)$ formé des $N-d$ premières colonnes de $\Delta(x)$ est de valuation $e''\leq c_\!Xe$ et que les autres mineurs $N-d\times N-d$ de cette matrice sont de valuation supérieure  (\textit{cf} lemme \ref{2}). Il nous suffit donc de prouver le résultat pour $\Delta_{e,e'}\cap h^{-1}(C)$. Soit $x\in \Delta_{e,e'}\cap h^{-1}(C)$. Soit $\varphi:\Spf R_{k'}\rightarrow\!X$ le $\mathbb{D}$-morphisme de schémas formels correspondant à $x$. Supposons $n\geq 2e$ et $n\geq e+c_\!Xe'$. Quitte à étendre les scalaires, on peut supposer que $k=\kappa(x)^{\mathrm{alg}}$. Nous allons procéder comme dans la preuve du lemme \ref{2}.

\textit{Commen\c{c}ons par étudier la fibre de $h_n$ au-dessus de $\pi_{n,\!X}(x)$}. Comme $\!Y$ est supposé lisse, on peut supposer qu'il existe un $\mathbb{D}$-morphisme étale de schémas formels $\!Y\rightarrow\mathbb{B}^d_R$ et, gr\^ace au lemme \ref{morphisme etale} que $\!Y=\mathbb{B}^d_R$ (\textit{cf} preuve du lemme \ref{4}). Par ailleurs l'assertion 1 et l'étude faite dans la preuve du lemme \ref{4} assurent que l'on peut également supposer que $\!X=\mathbb{B}^d_R$. Dans ce cas, cet ensemble peut se décrire comme
$$
\{\varphi + \pi^{n+1-e}u\mod \pi^{n+1} \mid (\!J_h(\varphi).u)\equiv 0\mod \pi^{e}\}.
$$
Comme $\!J_h(\varphi)$ est équivalente à une matrice diagonale, dont les termes diagonaux sont $\pi^{e_1},\ldots,\pi^{e_d}$ avec $\sum e_j=e$, cette description implique que la fibre de $h_n$ au-dessus de $\pi_{n,\!X}(h(\varphi))$ est $k$-isomorphe à $\mathbf{A}^{e}_k$.

On conclut alors par \textit{l'argument ``noethérien''} de la preuve du lemme \ref{2}.
\end{proof}

\begin{thm}
\label{tcv}
Soient $\!X$ et $\!Y$ deux $R$-schémas formels admissibles et réduits. Supposons que $\!Y$ est lisse. Soient $A$ et $B$ des ensembles fortement mesurables de $\G(\!X)$ et $\G(\!Y)$ respectivement. Supposons que $h$ induit une bijection entre $B$ et $A$ et que $h$ est tempéré sur $B$. Alors, pour toute application exponentiellement int\'egrable $\alpha:A\rightarrow\mathbf{Z}\cup\{\infty\}$, l'application $B\rightarrow\mathbf{Z}\cup\{\infty\}$ qui \`a $y$ associe $\alpha(h(y))+\ord_\pi(\mathrm{Jac})(y)$ est exponentiellement int\'egrable et on a la formule:
$$
\int_A\mathbf{L}^{-\alpha}d\mu=\int_B\mathbf{L}^{-\alpha\circ h-\ord_\pi(\mathrm{Jac})_h}d\mu.
$$
\end{thm}

\begin{proof}
Soit $\varepsilon>0$. Soit $A_I$ (\textit{resp.} $B_J$) une $\varepsilon$-approximation cylindrique de $A$ (\textit{resp.} $B$). Comme $h$ est tempéré sur $B$, on peut supposer que $B\cap\Sigma_h=\emptyset$. Prouvons tout d'abord que l'application $\beta:=\alpha\circ h+\ord_\pi(\mathrm{Jac})_h$ est mesurable sur $B$. Soit $n\in\mathbf{Z}$. Soit $F^{(n)}_L$ une $\varepsilon$-approximation cylindrique de $\alpha^{-1}(n)\subset A$. On va montrer que $(\alpha\circ h)^{-1}(n)$ est un ensemble mesurable de $\G(\!Y)$. On peut supposer que $F^{(n)}_0(\varepsilon)\subset A_0(\varepsilon)\subset A$ (\textit{cf} lemme \ref{approximation contenue}). On a l'égalité
$$
h^{-1}(\alpha^{-1}(n))\triangle h^{-1}(F^{(n)}_0(\varepsilon))\subset \cup_{i\in L} h^{-1}(F^{(n)}_i(\varepsilon)).
$$
Grâce au corollaire \ref{bidouille}, on peut supposer que $h^{-1}(A_0(\varepsilon))= B_0(\varepsilon)$. En particulier, cette remarque entra\^\i ne que
$$
h^{-1}(\alpha^{-1}(n))\triangle h^{-1}(F^{(n)}_0(\varepsilon))\subset (\cup_{j\in J} B_i(\varepsilon))\cup(\cup_{i\in L}h^{-1}(F^{(n)}_i(\varepsilon))\cap B_0(\varepsilon)).
$$
En effet, par hypothèse, on a 
$$
h^{-1}(\alpha^{-1}(n))\backslash B_0(\varepsilon)\subset B\cap B_0(\varepsilon),
$$
$$
(h^{-1}(\alpha^{-1}(n))\backslash h^{-1}(F^{(n)}_0(\varepsilon)))\cap B_0(\varepsilon)\subset \cup_{i\in L}h^{-1}(F^{(n)}_i(\varepsilon))\cap B_0(\varepsilon),
$$
et, comme $h^{-1}(F^{(n)}_0(\varepsilon))\subset B_0(\varepsilon)$,
$$
h^{-1}(F^{(n)}_0(\varepsilon))\backslash h^{-1}(\alpha^{-1}(n))\subset \cup_{i\in L}h^{-1}(F^{(n)}_i(\varepsilon))\cap B_0(\varepsilon).
$$
Le lemme \ref{5}/2 assure alors que $(\alpha\circ h)^{-1}(n)$ est mesurable. On a 
$$
\beta^{-1}(n)=\sqcup_{e\in\mathbf{N}} ((\alpha\circ h)^{-1}(n-e)\cap\ord_\pi(\mathrm{Jac})_h^{-1}(e)).
$$
La suite $(\mu_\!Y(\ord_\pi(\mathrm{Jac})_h^{-1}(e))_{e\in\mathbf{N}}$ tend vers 0. En effet, la suite $(\mu_\!Y(\pi_{e,\!Y}^{-1}(\G_e(\!Z_h))))_{e\in\mathbf{N}}$ tend vers $\mu_\!Y(\G(\!Z_h))$ (\textit{cf} proposition \ref{description des reunions de mesurables}/2); d'où la convergence vers 0 de la suite de terme général $\mu_\!Y(\ord_\pi(\mathrm{Jac})_h^{-1}(e))$. Le lemme \ref{mesurable comme limite} assure alors que l'application $\alpha\circ h +\ord_\pi(\mathrm{Jac})_h$ est mesurable sur $B$.

Il s'agit de montrer maintenant que l'intégrale converge et qu'on a bien la formule annoncée. Le lemme \ref{quasi-compacite de la topologie constructible} assure qu'il existe un entier $m\in\mathbf{N}$ tel que
$$
B_0(\varepsilon)\subset (\cup_{0\leq e'\leq m}h^{-1}(\G^{(e')}(\!X)))\cap(\cup_{0\leq e\leq m}\{y\in\G(\!Y)\mid\ord_\pi(\mathrm{Jac})_h(y)=e\}),
$$ 
Par additivit\'e de la mesure, il nous suffit de montrer le r\'esultat pour $B_0(\varepsilon)\subset\Delta_{e,e'}$ pour un entier $e\in\mathbf{N}$ et un entier $e'\in\mathbf{N}$. Soit $\varepsilon>0$. Soit $n\in\mathbf{Z}$. Notons $G_M^{(n-e)}$ une $\varepsilon$-approximation cylindrique de l'ensemble $(\alpha\circ h)^{-1}(n-e)$, qui est  mesurable dans $\G(\!X)$, telle que $G_0^{(n-e)}(\varepsilon)\subset B_0(\varepsilon)$. Il existe une $\varepsilon$-approximation cylindrique $F_L^{(n-e)}$ de $\alpha^{-1}(n-e)$ telle que $F_0^{(n-e)}(\varepsilon)=h(G_0^{(n-e)}(\varepsilon))$. Le lemme \ref{5} assure alors que:
$$
[\pi_{k,\!Y}( G_0^{(n-e)}(\varepsilon))]\mathbf{L}^{-n}=[\pi_{k,\!X}(F_0^{(n-e)}(\varepsilon))]\mathbf{L}^{-n+e}.
$$
La formule annoncée en découle.
\end{proof}

Grâce à la théorie des modèles de Néron faibles, on peut généraliser ce théorème au cas où seul la fibre générique de $\!Y$ est supposée lisse.

\begin{cor}
\label{tcvg}
Soient $\!X$ et $\!Y$ deux $R$-schémas formels admissibles et réduits. Supposons que $\!Y$ soit génériquement lisse. Soient $B$ et $A$ des ensembles fortement mesurables de $\G(\!Y)$ et $\G(\!X)$ respectivement, tels que $h$ induit une bijection de $B$ sur $A$. Supposons en outre que le morphisme $h$ soit tempéré sur $B$. Soit $\alpha: A\rightarrow \mathbf{Z}\cup\{\infty\}$ une application exponentiellement int\'egrable. Alors l'application qui \`a $y\in B$ associe $\alpha\circ h(y)+\ord_\pi(\mathrm{Jac})_h(y)$ est exponentiellement int\'egrable et l'on a l'\'egalit\'e dans $\widehat{\!M}$:
$$
\int_A\mathbf{L}^{-\alpha}d\mu_\!X=\int_B\mathbf{L}^{-(\alpha\circ h+\ord_\pi(\mathrm{Jac})_h)}d\mu_\!Y.
$$
\end{cor}

\begin{proof}
On peut supposer comme dans la preuve de \ref{tcv} que $B\cap \Delta_h=\emptyset$ et que $A\cap\G(\!X_{\mathrm{sing}})=\emptyset$. La théorie des modèles de Néron (\textit{cf} \cite{blr}, \cite{bs}, \cite{ls} \S 2.6) assure qu'il existe un ouvert $\!U$ quasi-compact d'un $R$-mod\`ele de la fibre g\'en\'erique de $\!Y$, not\'ee $Y_K$, lisse sur $R$, et un morphisme $g:\!U\rightarrow \!Y$ tel que $g$ induit une bijection entre les points de $\G(\!U)$ et ceux de $\G(\!Y)$. On a alors un diagramme commutatif:
$$
\xymatrix{\!U\ar[d]_{g}\ar@{-->}[dr]^{f}&\\
\!Y\ar[r]_{h}&\!X}
$$
Le morphisme $g$ induit une bijection entre $g^{-1}(B)$ et $B$ et \ref{hiv pour tcvg} assure que $g^{-1}(B)=f^{-1}(A)$ est fortement mesurable. Comme $g$ est tempéré, on peut supposer que $g^{-1}(B)(=f^{-1}(A))\cap \Sigma_g=\emptyset$. Il en résulte tout d'abord que 
$$
f^{-1}(A)\cap g^{-1}(h^{-1}(\G(\!X_{\mathrm{sing}})))=\emptyset.
$$ 
Par ailleurs, la formule de composition (\textit{cf} lemme \ref{formule de composition}) assure alors que $f^{-1}(A)\cap \Pi_f=\emptyset$. On a donc $f^{-1}(A)\cap \Sigma_f=\emptyset$ (en particulier, $f$ est tempéré sur $f^{-1}(A)$).

Supposons maintenant que l'application $\alpha\circ h +\ord_\pi(\mathrm{Jac})_h$ est exponentiellement int\'egrable sur $B$. Le th\'eor\`eme \ref{tcv} implique alors les \'egalit\'es dans $\widehat{\!M}$:
$$
\int_A\mathbf{L}^{-\alpha}d\mu_\!X=\int_{f^{-1}(A)}\mathbf{L}^{-\alpha\circ h\circ g-\ord_\pi(\mathrm{Jac})_{h\circ g}}d\mu_\!U
$$ 
et 
$$
\int_B\mathbf{L}^{-\alpha\circ h-\ord_\pi(\mathrm{Jac})_h}d\mu_\!Y=\int_{g^{-1}(B)}\mathbf{L}^{-\alpha\circ h\circ g-\ord_\pi(\mathrm{Jac})_h\circ g-\ord_\pi(\mathrm{Jac})_g}d\mu_\!U.
$$ 
La formule découle alors du lemme \ref{formule de composition}. Il faut prouver que l'application $\alpha\circ h +\ord_\pi(\mathrm{Jac})_h$ est exponentiellement intégrable sur $B$. Soit $F_n=(\alpha\circ h +\ord_\pi(\mathrm{Jac})_h)^{-1}(n)$. On remarque que 
$$
F_n=g(\sqcup_{e\in\mathbf{N}}(\alpha\circ h\circ g)^{-1}(n-e)\cap(\ord_\pi(\mathrm{Jac})_h\circ g)^{-1}(e)).
$$
Or la formule de composition \ref{formule de composition} assure que
$$
(\ord_\pi(\mathrm{Jac})_h\circ g)^{-1}(e)=\sqcup_{\mathrm{finie}}((\ord_\pi(\mathrm{Jac})_f)^{-1}(e-e')\cap(\ord_\pi(\mathrm{Jac})_g)^{-1}(e')).
$$
La mesurabilité de $F_n$ se déduit alors de la mesurabilité de $\alpha$ et du théorème \ref{theoreme de image directe}, comme dans la preuve de \ref{tcv}. La convergence de l'intégrale découle alors du lemme \ref{5}, comme dans la preuve du théorème \ref{tcv}.
\end{proof}

\bibliographystyle{amsplain}

\end{document}